\documentclass[12pt]{article}

\usepackage{amsthm,amssymb,mathtools} 
\usepackage{color}
\usepackage{enumitem}
\usepackage{hyperref}
\usepackage[subrefformat=parens]{subcaption}
\usepackage{mathrsfs} 

\usepackage{bm}
\usepackage{array}
\usepackage{arydshln}
\usepackage{comment}
\usepackage{graphicx}
\usepackage{mathrsfs}
\usepackage{here}
\usepackage{optidef}
\usepackage{url}

\usepackage{macros}
\usepackage{mathoperators}
\usepackage{theorems}

\usepackage{algorithmic}
\usepackage{algorithm}

\usepackage[capitalise]{cleveref}
\crefname{figure}{Figure}{Figures}
\crefname{table}{Table}{Tables}

\usepackage{tikz-cd} 

\hypersetup{
     colorlinks=true,
     linkcolor=black,
     filecolor=blue,
     citecolor = black,
     urlcolor=cyan,
     }

\def\0{\mbox{\bf 0}}
\def\1{\mbox{\bf 1}}
\def\2{\mbox{\bf 2}}
\def\3{\mbox{\bf 3}}
\def\4{\mbox{\bf 4}}
\def\5{\mbox{\bf 5}}
\def\6{\mbox{\bf 6}}
\def\7{\mbox{\bf 7}}
\def\8{\mbox{\bf 8}}
\def\9{\mbox{\bf 9}}

\def\c{\mbox{\boldmath $c$}}

\def\i{\mbox{\boldmath $i$}}

\def\l{\mbox{\boldmath $l$}}

\def\u{\mbox{\boldmath $u$}}

\def\AC{\mbox{$\cal A$}}
\def\BC{\mbox{$\cal B$}}
\def\CC{\mbox{$\cal C$}}
\def\DC{\mbox{$\cal D$}}

\def\IC{\mbox{$\cal I$}}

\def\MC{\mbox{$\cal M$}}

\def\OC{\mbox{$\cal O$}}

\def\SC{\mbox{$\cal S$}}
\def\TC{\mbox{$\cal T$}}

\def\XC{\mbox{$\cal X$}}

\def\ZC{\mbox{$\cal Z$}}

\def\Real{\mbox{$\mathbb{R}$}}
\def\Complex{\mbox{$\mathbb{C}$}}

\def\Natural{\mbox{$\mathbb{N}$}}

\def\bcirc{\mbox{$\mathrm{bcirc}$}}
\def\unfold{\mbox{$\mathrm{unfold}$}}
\def\fold{\mbox{$\mathrm{fold}$}}

\def\@themcountersep{}
\usepackage{fancybox,color}
\definecolor{lred}{rgb}{1,0.8,0.5}
\definecolor{lblue}{rgb}{0.8,0.8,1}
\definecolor{dred}{rgb}{0.6,0,0}
\definecolor{dblue}{rgb}{0,0,0.7}
\definecolor{violet}{rgb}{0.5804,0.0000,0.8275}
\definecolor{purple}{rgb}{0.2400,0.5700,0.2500}
\definecolor{TGreen}{rgb}{0,0.50,0.10}
\newcommand{\blue}[1]{\begin{color}{blue}#1\end{color}}

\usepackage[mathlines]{lineno} 
\usepackage{amsmath}

\usepackage{etoolbox} 

\usepackage{easy-todo}

\newcommand*\linenomathpatch[1]{%
    \cspreto{#1}{\linenomath}%
    \cspreto{#1*}{\linenomath}%
    \csappto{end#1}{\endlinenomath}%
    \csappto{end#1*}{\endlinenomath}%
}
\newcommand*\linenomathpatchAMS[1]{%
    \cspreto{#1}{\linenomathAMS}%
    \cspreto{#1*}{\linenomathAMS}%
    \csappto{end#1}{\endlinenomath}%
    \csappto{end#1*}{\endlinenomath}%
}

\expandafter\ifx\linenomath\linenomathWithnumbers
\let\linenomathAMS\linenomathWithnumbers
\patchcmd\linenomathAMS{\advance\postdisplaypenalty\linenopenalty}{}{}{}
\else
\let\linenomathAMS\linenomathNonumbers
\fi

\linenomathpatch{equation}
\linenomathpatchAMS{gather}
\linenomathpatchAMS{multline}
\linenomathpatchAMS{align}
\linenomathpatchAMS{alignat}
\linenomathpatchAMS{flalign}

\title{T-semidefinite programming relaxation with third-order tensors  for constrained polynomial optimization  } 

\makeatletter
\let\@fnsymbol\@arabic
\makeatother

\author{
\normalsize
    Hiroki Marumo\thanks{Department of  Mathematical and Computing Science,
        Tokyo Institute of Technology, 2-12-1-W8-29 Oh-Okayama, Meguro-ku, Tokyo 152-8552, Japan.
        ({\tt marumo.h.ac@m.titiech.ac.jp})
        ({\tt Makoto.Yamashita@c.titech.ac.jp}).
        The research of Makoto Yamashita was partially supported by JSPS KAKENHI Grant Number 21K11767 and 24K14836.}
\and
\normalsize
	Sunyoung Kim\thanks{Department of Mathematics, Ewha W. University, 52 Ewhayeodae-gil, Sudaemoon-gu,
	Seoul 	03760, Korea  ({\tt skim@ewha.ac.kr}). This work was supported
        by  NRF 2021-R1A2C1003810.}
\and
\normalsize
        Makoto Yamashita\footnotemark[1]
        }

\begin{document}
\maketitle


\begin{abstract} \noindent
We study T-semidefinite programming (SDP) relaxation for constrained polynomial optimization problems (POPs).
T-SDP relaxation for unconstrained POPs was introduced by Zheng, Huang and Hu in 2022.
In this work, we propose a T-SDP  relaxation for POPs with polynomial inequality constraints and show that the resulting T-SDP relaxation
formulated with third-order tensors can be transformed into the standard SDP relaxation with block-diagonal structures.
The convergence of the T-SDP relaxation 
to the optimal value of a given constrained POP
 is established under moderate assumptions as the relaxation level increases. 
Additionally, the feasibility and optimality of the T-SDP relaxation are  discussed.
Numerical results illustrate that the proposed T-SDP  relaxation enhances numerical efficiency.
\end{abstract}

\vspace{0.5cm}

\vspace{0.5cm}


\section{Introduction} \label{intro}

We consider constrained polynomial optimization problems (POPs):
        \begin{mini}
            {}
            {f(x)}
            {\tag{$\mathbb{P}_K$} \label{p_k}}
            {}
            \addConstraint{g_i(x)}{\ge 0,}{\quad i=1, \dots, r,}
        \end{mini}
        where $ f : \Real^n \to \Real $ and  $g_1, \dots, g_r : \Real^n \to \Real$  are real-valued polynomials.
Without loss of generality, we assume that no constant exists in the objective function.
The feasible set of \eqref{p_k} is denoted as $K := \left\{ x \in \Real^n \; \middle| \; g_i(x) \ge 0, \; \; i=1, \dots, r \right\}$.

Constrained POP \eqref{p_k} has been widely studied in the field of optimization
    such as 0-1 integer programming and non-convex quadratic optimization, and has various applications in medical
    bioengineering~\cite{ghosh2008polynomial, zhang2011new}, signal processing~\cite{mariere2003blind, qi2003multivariate},
    materials engineering~\cite{soare2008use}, and computer vision~\cite{hartley2003multiple}.

  As POP \eqref{p_k}  is known to be NP-hard,
 finding  a global optimal solution of \eqref{p_k}  is a significant challenge.
   A widely adopted strategy in solving  \eqref{p_k} has been the use of  semidefinite programming (SDP) relaxation based on Putinar's Positivstellensatz~\cite{putinar1993positive} by Lasserre~\cite{lasserre2001global}. Parrilo also proposed the SDP relaxation problem through the sum-of-squares (SOS) polynomial relaxation~\cite{parrilo2003semidefinite}, which is viewed as the dual of Lasserre's SDP relaxation.
    While Nesterov and Nemirovskii~\cite{nesterov1994interior} theoretically  showed that SDP problems can be solved in polynomial time, the
    size of SDP relaxation problems for POPs becomes increasingly large~\cite{nie2012regularization} 
    with the degree and the number of variables of POPs. Thus solving the SDP problems is challenging 
    due to computational costs~\cite{sun2020sdpnal+}.
    In fact,  for POPs of degree $d$ with $n$ variables, the size of the variable matrix in the SDP problem
     amounts to $s(d) :=  \binom{n+d}{d}$.
    To address this challenge, various approaches have been proposed, such as exploiting the sparse structures of  \eqref{p_k}~\cite{waki2006sums,WANG2021,QU2024} and
    employing weaker relaxations as diagonally-dominant-sum-of-squares relaxation~\cite{ahmadi2014dsos}. 

    For the difficulty in solving large-sized SDP relaxations of constrained POPs,
    our objective is to employ a third-order tensor and
    propose a relaxation method for constrained POP  \eqref{p_k}. 
    The third-order tensor has been studied in many applications~\cite{nion2009adaptive, tao2005stable}.
    Furthermore, Kilmer, Martin and Perrone~\cite{kilmer2008third} introduced the product in third-order tensors, an extension of the product of matrices, and third-order tensors have become gradually popular in recent years~\cite{du2022robust, martin2013order, semerci2014tensor}.

        Applying the product in third-order tensors, Zheng et al.~\cite{zheng2021t} introduced semidefiniteness and
         SDP in third-order tensors (referred to as semidefinite tensor and T-SDP, respectively).  
         Semidefinite tensor and T-SDP serve as extensions of   semidefinite matrices and SDPs, respectively. The entire set of semidefinite tensors is nonempty, closed, convex, pointed cone, and self-dual. Under these properties, duality theorems and optimality conditions have been shown for T-SDP as well as for SDP. While various application
         problems can be formulated as T-SDP,  no solver has currently been developed
         to directly solve the T-SDP problem. Therefore, it is common to use the structural advantage of third-order tensors to convert the problem into an equivalent SDP problem and solve it by applying  existing solvers such as \cite{aps2020mosek, yamashita2003implementation,yamashita2003sdpara}.

    As the close relationship between SOS relaxation and SOS polynomials,
     semidefinite third-order tensor and T-SDP are based on
     the SOS polynomial with a block circulant structure, commonly referred to as block circulant SOS polynomials~\cite{zheng2022unconstrained}. This block circulant SOS polynomial can be viewed as an extension of SOS polynomials.

The purpose of this paper is to introduce a relaxation by the block circulant SOS polynomial instead of the  SOS relaxation for constrained POPs
in \cite{kim2005generalized} 
to increase the efficiency of solving SDP relaxations of  \eqref{p_k}.
The SOS relaxation proposed in \cite{kim2005generalized} is called  the SOS relaxation or the basic SOS relaxation  in this paper.
Let the degree of $f(x)$ and $g_i(x)$  be $d$ and $w_i$, respectively, and 
          a positive integer $N$,  the relaxation level, such that $2N \ge d$ and $2N \ge \max_i \{w_i\}$,
          and  $\tilde{w}_i = \lceil w_i/2 \rceil,  \ i=1, \dots, r$, where $\lceil w \rceil$ denotes the smallest integer greater than equal to $w$.
With positive integers $m_0, l_0$, $m_i,$ and $ l_i$ such that
           $ s(N) = m_0 l_0, \ s(N-\tilde{w}_i) = m_i l_i, \ i=1, \dots, r$,
we employ $2N$-degree $l_0$-block circulant SOS polynomials 
	and $2(N-\tilde{w_i})$-degree
        $l_i$-block circulant SOS polynomials for our block-circulant SOS relaxation.
If we choose $ l_i = 1,\ i=0,\ldots,r$,
the proposed  block circulant SOS relaxation coincides with the basic SOS relaxation.
The proposed  block circulant SOS relaxation relies on the existence of
  $2N$-degree $l_0$ and $2(N-\tilde{w}_i)$-degree  $l_i$-block circulant SOS polynomials for \eqref{p_k}.
  Thus, the choices of $m_i \ \text{and} \ l_i,  i=0, \dots, r$ are essential for numerical efficiency, as shown in Section \ref{result}.

We first show that the block circulant SOS relaxation problem can be described as an equivalent T-SDP problem.
Then, the T-SDP relaxation problem is equivalently transformed into an SDP problem of smaller size than the SDP problem  from the basic SOS relaxation.
  Therefore,  the T-SDP relaxation is expected to reduce the computational time.
  Indeed, we observe the computational advantage through numerical experiments presented
in Section \ref{result}.
For instance,
the computational time for a constrained POP is reduced from 1264 seconds by the
basic SOS relaxation to 633 seconds by the proposed
relaxation based on circulant SOS polynomials with appropriate block sizes, under a mild assumption. 

The main contributions of this paper are:
\begin{itemize}

\item We propose a T-SDP relaxation for constrained POPs, extending the result in \cite{zheng2021t} for unconstrained POPs.
\item By transforming the T-SDP relaxation problem derived from a block circulant SOS relaxation for a constrained POP
 into the equivalent SDP problem of smaller size with the block-diagonal structure,
the computational efficiency of solving the constrained POP is increased. 

\item We analyze the feasibility and global optimality of the T-SDP relaxation problem, and prove the convergence of the T-SDP relaxation sequence to the optimal value of  \eqref{p_k} under an assumption that extends Putinar's Positivstellensatz.
\end{itemize}

    The remaining of our paper is organized as follows. In Section \ref{prel}, we describe the basic definitions of third-order tensor and block diagonalization, which are important for the equivalent transformation from T-SDP to SDP. Based on these definitions, we introduce
    the concepts of 
    semidefiniteness and  the T-SDP problem over third-order tensors and show that the T-SDP problem can be transformed into an equivalent SDP.
     We also define the block circulant SOS polynomial, an extension of the SOS polynomial, and introduce an important theorem that holds for semidefinite tensors.
    In Section \ref{prop}, we propose a block circulant SOS relaxation method that extends the basic SOS relaxation method for constrained POPs, and show that the relaxation problem can be described by a T-SDP form. We also show that it has an advantage in problem size compared to
    the  SDP relaxation. In addition, we provide theoretical analysis such as feasibility and global optimality for the T-SDP relaxation problem.
    In Section \ref{expe}, we compare the proposed T-SDP relaxation method  with the SDP relaxation using a set of
     test problems for constrained POPs.
    Finally, we conclude in Section \ref{conc}. 
    \vskip\baselineskip


    \section{Preliminaries}

        \label{prel}
        \subsection{Notation}
        \label{nota}
        Throughout this paper, $\Real^n$ denotes the set consisting of all real vectors of size $n$, $\Complex^n$  the set consisting of all complex vectors of size $n$, and $ \Natural^n$ the set of nonnegative integers.
        We use $\Real^{m \times n}$ to denote the set of $m \times n$ real matrices, $\Complex^{m \times n}$  the set of $m \times n$ complex matrices,
         $\Real^{m \times n \times l}$ the set of $m \times n \times l$ real third-order tensors. A third-order tensor is an array with one dimension added to the matrix.  
        We consider only real third-order tensors 
        in this paper.
        The $n$th-order identity matrix is denoted by $I_n$ and the $n$th-order zero matrix by $O_n$, or $I$ and $O$.
        The zero tensor
        is denoted by $\OC$. The identity tensor is defined in Section \ref{tensor}.

        The transpose of a matrix $A$ is denoted by $A^\mathsf{T}$, the conjugate by $\overline{A}$, and the conjugate transpose by $A^\mathsf{H}$.
        We denote the set consisting of all real symmetric (complex Hermitian) matrices of size $n \times n$ as  $\SymMat^{n \times n}$ ($\HerMat^{n \times n}$). 
        $A \succ(\succeq)   O$ denotes that 
        a matrix $A$ is positive definite (positive semidefinite), 
         while the set  of all real symmetric positive definite (positive semidefinite) matrices of size $n \times n$ is denoted by $\SymMat_{++}^{n \times n} (\SymMat_+^{n \times n})$. Similarly, we denote the set of all complex Hermitian positive definite (positive semidefinite) matrices of size $n \times n$ as $\HerMat_{++}^{n \times n}$ ($\HerMat_+^{n \times n})$.

        The inner product of two matrices in $\Real^{m \times n}$ is defined as $\langle A, B \rangle := \mathrm{Tr}(A^\mathsf{T}B) $ 
        $ = \sum_{i=1}^{m} \sum_{j=1}^{n}$ $ [A]_{(i,j)} [B]_{(i,j)}$, where $[A]_{(i,j)}$, $[B]_{(i,j)}$ denote the $(i, j)$th elements of $A$, $B$, respectively. The inner product of two third-order tensors $\AC$ and $\BC$ in $\Real^{m \times n \times l}$ is defined as 
        $\langle \AC, \BC \rangle := $ $ \sum_{i=1}^{m} \sum_{j=1}^{n} $ $ \sum_{k=1}^{l} [\AC]_{(i,j,k)} [\BC]_{(i,j,k)}$, where $[\AC]_{(i,j,k)}$, $[\BC]_{(i,j,k)}$ denote the $(i, j, k)$th elements of $\AC$, $\BC$, respectively. A symbol $\otimes$ denotes the Kronecker product of two matrices.

        \subsection{Third-order tensor}
        \label{tensor}
       For a third-order tensor $\AC \in \Real^{m \times n \times l}$, we consider $\AC$ as a stack of frontal slices $A^{(i)} \in \Real^{m \times n}, i=1, \dots, l$, as proposed in \cite{kilmer2008third} where several operators on
         $\AC  \in \Real^{m \times n \times l}$ were introduced as follows:
        \begin{align*}
            &\bcirc(\AC) :=
            \left[ \begin{array}{@{\,}ccccc@{\,}}
                A^{(1)} & A^{(l)}   & A^{(l-1)} & \ldots & A^{(2)} \\
                A^{(2)} & A^{(1)}   & A^{(l)}   & \ldots & A^{(3)} \\
                A^{(3)} & A^{(2)}   & A^{(1)}   & \ldots & A^{(4)} \\
                \vdots  & \vdots    & \vdots    & \ddots & \vdots \\
                A^{(l)} & A^{(l-1)} & A^{(l-2)} & \ldots & A^{(1)}
            \end{array} \right], \quad
            \unfold(\AC) :=
            \left[ \begin{array}{@{\,}c@{\,}}
                A^{(1)} \\
                A^{(2)} \\
                A^{(3)} \\
                \vdots \\
                A^{(l)}
            \end{array} \right], \\
            &\bcirc^{-1}_{l}(\bcirc(\AC)) := \AC, \quad \fold_{l}(\unfold(\AC)) := \AC.
        \end{align*}
        Let $\bcirc : \Real^{m \times n \times l} \to \Real^{ml \times nl}$ be the operator that arranges each frontal slice of $\AC$ into a block circulant matrix, $\unfold : \Real^{m \times n \times l} \to \Real^{ml \times n}$ be the operator that arranges each frontal slice of $\AC$ in columns. $\bcirc^{-1}_{l}$ and $\fold_{l}$ are operators that are the inverse operations of $\bcirc$ and $\unfold$ with $l$, respectively. 
From the definition of $bcirc$, it clearly holds that
\begin{equation}
    \langle \AC, \BC \rangle = \frac{1}{l} \left\langle \bcirc(\AC), \bcirc(\BC) \right\rangle. \label{Lemma1}
\end{equation}

        It is known that any circulant matrix $A \in \Real^{n \times n}$ can be diagonalized with the normalized
        discrete Fourier transform (DFT) matrix, where the DFT matrix is the Fourier matrix of size $n \times n$ defined as
        \begin{equation*} F_n :=
            \frac{1}{\sqrt{n}}
        \left[ \begin{array}{@{\,}ccccc@{\,}}
            1      & 1            & 1               & \ldots & 1 \\
            1      & \omega       & \omega^2        & \ldots & \omega^{n-1} \\
            1      & \omega^2     & \omega^4        & \ldots & \omega^{2(n-1)} \\
            \vdots & \vdots       & \vdots          & \ddots & \vdots \\
            1      & \omega^{n-1} & \omega^{(n-1)2} & \ldots & \omega^{(n-1)(n-1)}
        \end{array} \right], \quad \omega := e^{\frac{2\pi i}{n}}.
        \end{equation*}
        Using this, the block circulant matrix  $\bcirc(\AC)$ can be block-diagonalized with the Kronecker product.
        
      \begin{lemma}  \label{Lemma2}\cite{kilmer2013third} Any third-order tensor $\AC \in \Real^{m \times n \times l}$ can be block-diagonalized as
        \begin{equation*}
            (F_l^\mathsf{H} \otimes I_m) \bcirc(\AC) (F_l \otimes I_n) = Diag(A_1, A_2, \dots, A_l) =
            \left[ \begin{array}{@{\,}cccc@{\,}}
                A_1 &     &        & \\
                    & A_2 &        & \\
                    &     & \ddots & \\
                    &     &        & A_l
            \end{array} \right],
        \end{equation*}
        where $A_1 \in \Real^{m \times n}$, $A_i \in \Complex^{m \times n}$ and $A_i = \overline{A_{l+2-i}}, \, i=2,\dots, l$.
        \end{lemma}

        In what follows, we define the product and transposition of third-order tensors and extend the
        semidefiniteness to the space of third-order tensors,
        and introduce some theorems on semidefinite tensors.
     
        \begin{definition}
        \cite[Definition 4.1]{kilmer2008third} Let $\AC \in \Real^{m \times n \times l}$ and $\BC \in \Real^{n \times p \times l}$ be two third-order tensors. Then, the product of third order tensors $\AC * \BC \in \Real^{m \times p \times l}$ is defined by
        \[ \AC * \BC := \fold_l\left( \bcirc(\AC) \unfold(\BC )\right). \]
        \end{definition}
        In the case of $l = 1$, the third-order tensors $\AC \in \Real^{m \times n \times l}$, $\BC \in \Real^{n \times p \times l}$ are  the matrices $A \in \Real^{m \times n}$, $B \in \Real^{n \times p}$, respectively, and
         $\bcirc(A) = A$ and $\unfold(B) = B$ hold for $l = 1$, thus the product in the third-order tensors $\AC * \BC$ is the matrix product $AB$.

                We call a third-order tensor whose  first frontal slice is the identity matrix and the other frontal slices are the zero matrix
        the identity tensor $\IC$.
        For a third-order tensor $\AC \in \Real^{n \times n \times l}$ and the
        identity tensor $\IC \in \Real^{n \times n \times l}$, we have $\AC * \IC = \IC * \AC = \AC$. 
    
        \begin{definition} \cite[Definition 4.7]{kilmer2008third} If $\AC \in \Real^{m \times n \times l}$ is a third-order tensor, then the transpose $\AC^\mathsf{T}$ is obtained by transposing each of the frontal slices $A^{(i)}$ and then reversing the order of transposed frontal slices $2$ through $l$. Furthermore, $\AC \in \Real^{m \times m \times l}$ is called a symmetric sensor if $\AC^\mathsf{T} = \AC$. We denote the set of all real symmetric tensors of size $m \times m \times l$ as $\SymMat^{m \times m \times l}$.
             \end{definition}

       \begin{definition}
        \cite[Definition 6]{zheng2021t} Let $\AC \in \SymMat^{m \times m \times l}$ be a symmetric tensor. We say $\AC$ is symmetric 
        positive (semi)definite, if and only if,
        \begin{align*}
            \langle \XC, \AC*\XC \rangle > (\ge) 0
        \end{align*}
        holds for any $\XC \in \Real^{m \times 1 \times l} \setminus \{\OC\}$(for any $\XC \in \Real^{m \times 1 \times l}$).
        \end{definition}
         We denote $\AC \succ_{\TC} (\succeq_{\TC}) \ \OC$  if $\AC \in \SymMat^{m \times m \times l}$ is positive (semi)definite, and 
         $\SymMat_{++}^{m \times m \times l} \ (\SymMat_+^{m \times m \times l})$  denotes the set  of all real symmetric positive
         (semi)definite tensors of size $m \times m \times l$.         
        The set $\SymMat_+^{m \times m \times l}$  is a nonempty, closed, convex, and pointed cone \cite[Proposition 3]{zheng2021t}. Moreover, 
        $\SymMat_+^{m \times m \times l} = (\SymMat_+^{m \times m \times l})^*$ holds for the dual cone $(\SymMat_+^{m \times m \times l})^*$ of $\SymMat_+^{m \times m \times l}$ (the self-duality) \cite[Theorem 9]{zheng2021t}.

              \begin{theorem} \label{Theorem1}
        \cite[Theorem 4]{zheng2021t} Let $\AC \in \SymMat^{m \times m \times l}$ be a symmetric tensor. Then, the following statements are equivalent: \vspace{-2mm}
        \begin{itemize}
            \item[(i)] $\AC \in \SymMat_+^{m \times m \times l}$.
            \item[(ii)] $\bcirc(\AC) \in \SymMat_+^{ml \times ml}$.
            \item[(iii)] The block-diagonal matrix $Diag(A_1, \dots, A_l)$, which is a block diagonalization of A, is a Hermitian positive semidefinite matrix. In other words, $A_1 \in \SymMat^{m \times m}_+$, $A_i \in \HerMat^{m \times m}_+, \, i=2, \dots, l$.
        \end{itemize}
        \end{theorem}
   
        \subsection{SDP in third-order tensor space}
        \label{t-sdp}
       We discuss 
        an SDP  problem in the space of third order tensors (T-SDP) in this section 
       using the semidefiniteness in third-order tensor defined in Section \ref{tensor}. We also describe that a T-SDP problem
        can be converted into
       an SDP problem in complex space.
       The transformed SDP can be solved using  existing SDP
       solvers, for instance Mosek \cite{aps2023mosek}, which is used to solve the  SDP  
       in our numerical experiments  in Section \ref{result}.
 
        An SDP problem is defined as: 
        \begin{align*}
        \mbox{ minimize } \left\langle C, X \right\rangle \ \mbox{subject to } \ \ \left\langle A_i, X \right\rangle =b_i, \quad i=1, \dots, r, \ \ X\succeq O,
        \end{align*}
        where $X \in \SymMat^{m \times m}$ is the decision variable, $C, A_i \in \SymMat^{m \times m}$ and $b \in \Real^{r}$ are given symmetric matrices 
        and a vector, respectively.
        Its dual problem is described as:
        \begin{align*}
          \mbox{ maximize } \sum_{i=1}^{r} b_iy_i \ \ \mbox{subject to } \ \ \sum_{i=1}^{r} y_i A_i + S =C, \ \ S\succeq O,
         \end{align*}
        where $y \in \Real^{r}$ is the decision variable and $S \in \SymMat^{m \times m}$ is the slack variable.

        Now, we  extend the SDP problem in the space of symmetric matrices
        to the space of third-order tensors by extending the real symmetric matrices $C, A_i, X$ to  real symmetric tensors $\CC, \AC_i, \XC \in \SymMat^{m \times m \times l}$ and the semidefinite constraint $X \succeq O$ to $\XC \succeq_{\TC} \OC$.
        A T-SDP problem is defined as: 
         \begin{align} \tag*{(PTSDP)} \label{ptsdp}
        \mbox{ minimize } \left\langle \CC, \XC \right\rangle \ \mbox{subject to } \ \ \left\langle \AC_i, \XC \right\rangle =b_i, \quad i=1, \dots, r, \ \ \XC\succeq_{\TC} \OC.
        \end{align}
        Similarly, its dual problem is defined as: 
        \begin{align} \tag*{(DTSDP)} \label{dtsdp}
          \mbox{ maximize } \sum_{i=1}^{r} b_iy_i \ \ \mbox{subject to } \ \ \sum_{i=1}^{r} y_i \AC_i + \SC= \CC, \ \ \SC \succeq_{\TC} \OC.
        \end{align}
        Clearly, an SDP problem can be considered as a special case of T-SDP  at $l=1$.

        The following duality theorem holds for T-SDP as well as SDP.
        \begin{theorem} \label{Theorem2}
        \cite[Theorem 11]{zheng2021t} Let $F(P)$, $F(D)$, $p^*$, and $d^*$ be defined as follows:
        \begin{align*}
            & F(P) = \left\{ \XC \in \SymMat^{m \times m \times l} \; \middle| \; \left\langle \AC_i,  \XC \right\rangle = b_i,\, i=1, \dots, r, \;\XC \succeq_{\TC} \OC \right\}, \\
            & F(D) = \left\{ (y, \SC) \in \Real^r \times \SymMat^{m \times m \times l} \; \middle| \; \sum_{i=1}^{r} y_i \AC_i + \SC= \CC,\; \SC \succeq_{\TC} \OC \right\}, \\
            & p^* = \inf \left \{ \langle \CC , \XC \rangle \mid \XC \in F(P) \right \}, \quad d^* = \sup \left \{ \langle b , y \rangle \mid (y, \SC) \in F(D) \right \}.
        \end{align*}
        Suppose that $\XC \in F(P)$ and $(y, \SC) \in F(D)$. Then, 
        $\langle b, y \rangle \le \langle \CC, \XC \rangle$. In addition,
        if  (PTSDP) is bounded below and strictly feasible or (DTSDP)  is bounded above and strictly feasible, 
        then the other  is solvable and  $p^*=d^*$.
        \end{theorem}
         Based on \cite{zheng2021t}, we briefly describe a method for solving the T-SDP problem by transforming it into an equivalent SDP
        problem over the complex space. 
        For every third-order tensor $\CC \in \SymMat^{m \times m \times l}$, $\bcirc(\CC)$ can be block-diagonalized as
        \begin{equation*}
            \bcirc(\CC) = (F_l \otimes I_m) Diag(C_1, C_2, \dots, C_l) (F_l^\mathsf{H} \otimes I_m).
        \end{equation*}
        Thus, for the objective function of \ref{ptsdp}, we obtain
        \begin{align*}
            \langle \CC, \XC \rangle & = \frac{1}{l} \langle \bcirc(\CC), \bcirc(\XC) \rangle 
             = \frac{1}{l} \mathrm{Tr}(\bcirc(\CC)\bcirc(\XC)) \\
            & = \frac{1}{l} \mathrm{Tr}( (F_l \otimes I_m) Diag(C_1, \dots, C_l) (F_l^\mathsf{H} \otimes I_m) (F_l \otimes I_m) Diag(X_1, \dots, X_l) (F_l^\mathsf{H} \otimes I_m) ) \\
            & = \frac{1}{l} \mathrm{Tr}(Diag(C_1, \dots, C_l) Diag(X_1, \dots, X_l)) 
             = \frac{1}{l} \sum_{k=1}^{l} \langle C_k, X_k \rangle,
        \end{align*}
        using \eqref{Lemma1}  for the first equality and the commutative property of trace for the fourth equality.

        Similarly, using the block-diagonalized matrix $Diag(A^i_1, \dots, A^i_l)$ of $\AC_i$, we have
         $\langle \AC_i, \XC \rangle = \frac{1}{l} \sum_{k=1}^{l} \langle A^i_k, X_k \rangle$. From Theorem \ref{Theorem1}, $\XC \succeq_{\TC} \OC$ is equivalent to $X_1, \dots, X_l \succeq O$.
        Since $l$ is an integer such that $l \ge 1$, we can replace $\frac{1}{l} X_k$ by $X_k$ without loss of generality. 
        Consequently, \ref{ptsdp} and \ref{dtsdp} are transformed into SDP problems in the space of complex matrices:
        \vspace{-2mm}
        \begin{align} \tag*{(PCSDP)} \label{pcsdp}
        \mbox{ min } \sum_{k=1}^{l} \langle C_k, X_k \rangle \ \mbox{sub. to } \ \ \sum_{k=1}^{l} \langle A^i_k, X_k \rangle=b_i, \quad i=1, \dots, r, \ \ 		X_1, \dots, X_l \succeq O,
        \end{align}
            \vspace{-2mm}
\vspace{-4mm}
        \begin{align} \tag*{(DCSDP)} \label{dcsdp}
          \mbox{ max } \sum_{i=1}^{r} b_iy_i \ \ \mbox{sub. to } \ \ \sum_{i=1}^{r} y_i A^i_k + S_k=C_k, \quad k=1, \dots, l, \ \ S_1, \dots, S_l \succeq O,
         \end{align}
        where $X_k \in \HerMat^{m \times m}$, $y \in \Real^{r}$ and $ S_k \in \HerMat^{m \times m}$ are the decision variables,
         and $C_k, $ $ A^i_k \in \HerMat^{m \times m}$
         are obtained from the diagonalization.

Since Lemma~\ref{Lemma2} holds for block diagonal matrices as described in Section \ref{tensor}, \ref{pcsdp} can be transformed into an equivalent problem of smaller size depending on the values of $l$, even or odd. \\ 
%
        \textbf{Case 1: $l$ is even} \\
        Let $C_k, A^i_k, X_k, \, k=1, \dots, l$ be the blocks in the block diagonalized matrix of $\CC, \AC_i, \XC$ in \ref{ptsdp}.
        Then,  the following relation holds for $C_k, A^i_k, X_k$ by \eqref{Lemma1}: 
        \begin{align*}
            &C_1 \in \SymMat^{m \times m}, \, C_{\frac{l+2}{2}} \in \SymMat^{m \times m}, \, C_k \in \HerMat^{m \times m}, \, C_k = \overline{C_{l+2-k}}, \\
            &A^i_1 \in \SymMat^{m \times m}, \, A^i_{\frac{l+2}{2}} \in \SymMat^{m \times m}, \, A^i_k \in \HerMat^{m \times m}, \, A^i_k = \overline{A^i_{l+2-k}}, \quad i=1. \dots, r, \\
            &X_1 \in \SymMat^{m \times m}, \, X_{\frac{l+2}{2}} \in \SymMat^{m \times m}, \, X_k \in \HerMat^{m \times m}, \, X_k = \overline{X_{l+2-k}}.
        \end{align*}
        For the objective function of \ref{pcsdp}, we obtain
        \begin{align*}
           \sum_{k=1}^{l} \left\langle C_k, X_k \right\rangle &= \left\langle C_1, X_1 \right\rangle + \left\langle C_{\frac{l+2}{2}}, X_{\frac{l+2}{2}} \right\rangle + \sum_{k=2}^{\frac{l}{2}} \left( \left\langle C_k, X_k \right\rangle + \left\langle \overline{C_k}, \overline{X_k} \right\rangle \right) \\
           &= \left\langle C_1, X_1 \right\rangle + \left\langle C_{\frac{l+2}{2}}, X_{\frac{l+2}{2}} \right\rangle + 2\sum_{k=2}^{\frac{l}{2}} \left\langle C_k, X_k \right\rangle,
        \end{align*}
        using the fact that the inner product of the two Hermitian matrices is equal to the inner product of their conjugates in the  second equality.
        Similarly, we obtain
        \begin{align*}
           \sum_{k=1}^{l} \left\langle A^i_k, X_k \right\rangle = \left\langle A^i_1, X_1 \right\rangle + \left\langle A^i_{\frac{l+2}{2}}, X_{\frac{l+2}{2}} \right\rangle + 2\sum_{k=2}^{\frac{l}{2}} \left\langle A^i_k, X_k \right\rangle, \quad i=1, \dots, r.
        \end{align*}
        Therefore, \ref{pcsdp} is equivalent to the following problem:
        \begin{mini}
            {}  
            {\left\langle C_1, X_1 \right\rangle + 2\sum_{k=2}^{\frac{l}{2}} \left\langle C_k, X_k \right\rangle + \left\langle C_{\frac{l+2}{2}}, X_{\frac{l+2}{2}} \right\rangle} 
            {\tag*{(P'CSDP)} \label{p'csdp}}  
            {}  
            \addConstraint{\left\langle A^i_1, X_1 \right\rangle + 2\sum_{k=2}^{\frac{l}{2}} \left\langle A^i_k, X_k \right\rangle + \left\langle A^i_{\frac{l+2}{2}}, X_{\frac{l+2}{2}} \right\rangle}{=b_i,} {\quad i=1, \dots, r} 
            \addConstraint{X_1, \dots, X_{\frac{l+2}{2}}}{\succeq O.} 
        \end{mini}
%
        \textbf{Case 2: $l$ is odd} \\
        As in Case 1, the following relation holds for $C_k, A^i_k, X_k$ by \eqref{Lemma1}:
        \begin{align*}
            &C_1 \in \SymMat^{m \times m}, \, C_k \in \HerMat^{m \times m}, \, C_k = \overline{C_{l+2-k}}, \\
            &A^i_1 \in \SymMat^{m \times m}, \, A^i_k \in \HerMat^{m \times m}, \, A^i_k = \overline{A^i_{l+2-k}}, \quad i=1. \dots, r, \\
            &X_1 \in \SymMat^{m \times m}, \, X_k \in  \HerMat^{n \times n}, \, X_k = \overline{X_{l+2-k}}.
        \end{align*}
        Therefore, \ref{pcsdp} is equivalent to the following problem:
        \begin{mini}
            {}  
            {\left\langle C_1, X_1 \right\rangle + 2\sum_{k=2}^{\frac{l+1}{2}} \left\langle C_k, X_k \right\rangle} 
            {\tag*{(P''CSDP)} \label{p''csdp}}  
            {}  
            \addConstraint{\left\langle A^i_1, X_1 \right\rangle + 2\sum_{k=2}^{\frac{l+1}{2}} \left\langle A^i_k, X_k \right\rangle}{=b_i,} {\quad i=1, \dots, r} 
            \addConstraint{X_1, \dots, X_{\frac{l+1}{2}}}{\succeq O.} 
        \end{mini}

        We have shown that \ref{ptsdp} with a third-order tensor of size $m \times m \times l$ can be transformed into \ref{pcsdp} which 
        includes the sum of $l$ third-order tensors of size $m \times m$ matrices. Furthermore, \ref{pcsdp} can be transformed into \ref{p'csdp} and \ref{p''csdp} which include the sum of $\frac{l+2}{2}$ or $\frac{l+1}{2}$ matrices of size $m \times m$, respectively, depending on even or odd value of $l$.
  These properties will be used to propose a T-SDP relaxation for constrained POPs.
        \subsection{Block circulant SOS polynomials}
        \label{bcircSOSpoly}
      In this section, we first define the SOS polynomial with block circulant structure,
      referred to as block circulant SOS polynomial, and then explain the relationship between the block circulant SOS polynomial and semidefinite tensor.

We first describe some notation.
A real-valued polynomial $f(x)$ of degree $d$, $d\in \Natural$, is expressed as
\begin{gather}
            f(x) = \sum_{0 \le |\alpha| \le d} b_{\alpha} x^{\alpha} \ \text{with} \
            x^{\alpha} := x_1^{\alpha_1} x_2^{\alpha_2} \cdots x_n^{\alpha_n}, \ b_\alpha := b_{\alpha_1 \alpha_2 \cdots \alpha_n} \in \Real \ \text{and} \; |\alpha| := \sum_{i=1}^{n} \alpha_i \le d,
            \label{express-f}
        \end{gather}
for $x \in \Real^n$ and $\alpha \in \Natural^n$.
        Let
        \begin{equation*}
            [x]_d := [1, x_1, \dots, x_n, x_1^2, x_1x_2, \dots, x_n^2, \dots, x_1^d, \dots, x_n^d ]^\mathsf{T},
        \end{equation*}
        and let $s(d) = \binom{n+d}{d}$ be the dimension of $ [x]_d$.
        Let $\Real[x]_d$ be the set of all $d$-degree real-valued polynomials.
        If $f(x)  \in \Real[x]_d $ and $f(x) \ge 0 $ for any $x \in \Real^n$, then
        $f $ is called a nonnegative polynomial. The set of all nonnegative polynomials is denoted by $\Real_+[x]_d$.
        We say that $q$ is an SOS polynomial if the real-valued polynomial $q : \Real^n \to \Real$ can be expressed as $q(x) = \sum_{i} \bar{q}_i^2(x)$ 
        with some real-valued polynomial $\bar{q}_i : \Real^n \to \Real$.
         It is well-known that a nonnegative polynomial is not necessarily an SOS polynomial except for special cases \cite{reznick2000some}.

        \begin{definition} \label{Definition4}
         \cite[Definition 3]{zheng2022unconstrained} Let $q: \Real^n \to \Real$ be a  real-valued polynomial
         of degree $2d$, and $m$ and $l$ be positive integers such that $s(d) = m \cdot l$. 
            Then we say $q$ is an $l$-block circulant SOS polynomial if 
        \begin{equation*}
            q(x) = \sum_{i=1}^{r} \sum_{j=1}^{l} \left( \left(q_j^i \right)^\mathsf{T} [x]_d \right)^2,
        \end{equation*}
        where $q_j^i$ is the $j$-th column vector of $\bcirc(\fold_l(q^i))$ for some generators $\{q^1, \dots, q^r\} \subset \Real^{s(d)}$.   
        \end{definition}
       \begin{example} \label{example:block-ciculant}
        ($3$-block circulant SOS polynomial) \\
     A polynomial  $q$ in $2$ variables ($x_1$ and $x_2$) defined by 
       \begin{align*}
           q(x) = 21 + 8 x_2 + 9x_1^2 + 20 x_1 x_2 + 33 x_2^2 + 8 x_1 x_2^2 + 9 x_1^4 + 12 x_1^3 x_2 + 21 x_1^2 x_2^2 + 9 x_2^4 
       \end{align*}
        is a $3$-block circulant SOS polynomial.
        Two generators $q^1, q^2 $ for $q(x)$ are 
       \begin{align*}
           q^1 = [2, 0, 0, 0, 0, 3]^\mathsf{T}, \quad q^2 = [4, 0, 1, 0, 0, 0]^\mathsf{T}.
       \end{align*}
       We see that
       \begin{align*}
           \bcirc(\fold_3(q^1)) = \left[ \begin{array}{@{\,}c:c:c@{\,}}
                   2 & 0 & 0 \\
                   0 & 3 & 0 \\ \hdashline
                   0 & 2 & 0 \\
                   0 & 0 & 3 \\ \hdashline
                   0 & 0 & 2 \\
                   3 & 0 & 0
               \end{array} \right], \quad
           \bcirc(\fold_3(q^2)) = \left[ \begin{array}{@{\,}c:c:c@{\,}}
                   4 & 0 & 1 \\
                   0 & 0 & 0 \\ \hdashline
                   1 & 4 & 0 \\
                   0 & 0 & 0 \\ \hdashline
                   0 & 1 & 4 \\
                   0 & 0 & 0
               \end{array} \right].
       \end{align*}
       From
       $ [x]_2 := [1, x_1, x_2, x_1^2, x_1x_2,  x_2^2 ]^\mathsf{T}$, 
       we have
       \begin{align*}
           \sum_{i=1}^{2} \sum_{j=1}^{3} ((q_j^i)^\mathsf{T} [x]_2)^2 = 
           & \ (2 + 3 x_2^2)^2 + (3 x_1 + 2 x_2)^2 + (3 x_1^2 + 2 x_1 x_2 )^2  \\
           & \ + (4 + x_2)^2 + (4 x_2 + x_1 x_2)^2 + (1 + 4 x_1 x_2)^2,
       \end{align*}
       which is equivalent to $q$.
       \end{example}

       If $l = 1$   in Definition \ref{Definition4}, then ${\rm bcirc}({ \rm fold}_l(q^i)) = q^i$, thus $q_j^i$ is $q^i$ itself.
        Thus, a 1-block circulant SOS polynomial is an SOS polynomial.


        \begin{theorem} \label{Theorem3}
        \cite{zheng2022unconstrained} Let $q: \Real^n \to \Real$ be a real-valued polynomial of  degree $2d$, and $m$ and $l$ be positive integers such that $s(d) = m \cdot l$. Then, the following statements are equivalent:
        \begin{itemize}
            \item[(i)] $q$ is an $l$-block circulant SOS polynomial
            \item[(ii)] There exists an $l$-block circulant matrix $A \in \SymMat^{ml \times ml}_+$ such that each block 
            $A^{(1)}, \dots, A^{(l)}\in \Real^{m \times m}$
             and  $q$ can be expressed by
            \begin{equation*}
                q(x) = \left\langle A, [x]_d
                [x]_d^\mathsf{T} \right\rangle
            \end{equation*}
            \item[(iii)] There exists a semidefinite tensor $\AC \in \SymMat^{m \times m \times l}_+$ and $[\XC]_d^l = \fold_l([x]_d)$ such that $q$ can be expressed by
            \begin{equation*}
                q(x) = \left\langle \AC, [\XC]_d^l *
                \left( [\XC]_d^l \right)^\mathsf{T} \right\rangle.
            \end{equation*}
        \end{itemize}
        \end{theorem}

        In the case of $l = 1$ in Theorem \ref{Theorem3}, the semidefinite $1$-block circulant matrix in $(ii)$ is  a semidefinite matrix. Moreover, since the third-order tensor $\AC$ in $(iii)$ becomes a matrix and $[\XC]_d^l$  a vector $[x]_d$, $(ii)$ and $(iii)$ are equivalent. 
        For $l \ge 2$, the $l$-block circulant matrix $A$ in $(ii)$ can be regarded as a $1$-block circulant matrix when the entire matrix is considered as a single block. As a result, an $l$-block circulant SOS polynomial is also an SOS polynomial. Consequently, we have the following relation:
       \begin{align}
            \Sigma^l[x]_d \, (l \ge 2) \subseteq \, \Sigma^1[x]_d \, \subseteq \, \Real_+[x]_d, \label{eq:Inclusion}
        \end{align}
        where $ \Sigma^l[x]_d$ denotes the set of $l$-block circulant SOS polynomials of degree $d$.

%
   \section{T-SDP relaxations for constrained POPs}
    \label{prop}
  In this section, we propose a T-SDP relaxation method for the constrained POP \eqref{p_k}, which is an extension of the  SDP relaxation
  method in Section~\ref{t-sdp}.
   Then,
   we show that the proposed T-SDP relaxation  can be reduced to a smaller SDP using the properties of T-SDP
     than the SDP problem from
     the basic SOS relaxation.
   We also discuss the feasibility and global optimality of the proposed method.

        \subsection{SDP relaxation in the third-order tensor space}
        \label{tsdp-relax}
	To derive T-SDP relaxation for \eqref{p_k}, we
        let $\tilde{w}_i = \lceil w_i/2 \rceil$, where $w_i$ is the degree of $g_i(x)$ and
         define a positive integer $N$ called the relaxation level such that $2N \ge d$ and $2N \ge \max_i \{w_i\}$.
         We also determine certain positive integers $m_0, l_0$, $m_i,$ and $ l_i, \,  \ i=1,\dots, r$ such that
        \begin{align}
            \label{mlsplitting}
            s(N) = m_0 l_0, \quad s(N-\tilde{w}_i) = m_i l_i, \quad i=1, \dots, r.
        \end{align}
  
        For the discussion of  the optimal value of \eqref{p_k}, we introduce the following assumptions.
        We should mention that the convexity and connectivity of the feasible set $K$ are not assumed here.

\begin{assum} \label{Assumption1}
        \cite[Assumption 4.1]{lasserre2001global} The set $K = \left\{ x \in \Real^n \; \middle| \; g_i(x) \ge 0, \; \; i=1, \dots, r \right\}$ is compact and
        there exists a real-valued polynomial $u : \Real^n \to \Real$ such that $\{ x \in \Real^n \mid u(x) \ge 0 \}$ is compact, then
        \begin{align}
            \label{assumption1}
            u(x) = u_0(x) + \sum_{i=1}^{r} g_i(x) u_i(x) \quad \text{for all } x \in \Real^n,
        \end{align}
        where $u_i(x), \, i=0, \dots, r$ are SOS polynomials.
\end{assum}
        \begin{assum} \label{Assumption2}
        The set $K = \left\{ x \in \Real^n \; \middle| \; g_i(x) \ge 0, \; \; i=1, \dots, r \right\}$ is compact. If $p^*_K = \min_{x \in K} f(x)$, for any $\epsilon > 0$,
        \begin{align}
            \label{assumption2}
            f(x) - p_K^* + \epsilon = q_0(x) + \sum_{i=1}^{r} g_i(x) q_i(x) \quad \text{for all } x \in K,
        \end{align}
        where $q_0(x)$ is a $2N$-degree $l_0$-block circulant SOS polynomial, $q_i(x), \, i=1, \dots, r$ are $2(N-\tilde{w}_i)$-degree
        $l_i$-block circulant SOS polynomials.
        \end{assum}
       We note that Assumption \ref{Assumption1} is commonly
        used
        and Assumption  \ref{Assumption2} is an extension of Putinar's Positivstellensatz \cite{putinar1993positive}, which was used in the SDP
        relaxation  \cite{lasserre2001global}
         for the objective function $f$. In Assumption~\ref{Assumption2}, $q_0$ and $q_i$ need to be $l_0$- and $l_i$-block circulant SOS polynomials
          instead of  SOS polynomials. 
         Thus, Assumption \ref{Assumption2} always holds if 
         $l_0 = 1$ and $l_i = 1$ under the condition
         for Putinar's Positivstellensatz. 

        We extend the basic SOS relaxation method based on the generalized Lagrange function by Kim et al.~\cite{kim2005generalized} to
         SOS relaxations using block circulant SOS polynomials.
        The generalized Lagrangian function for \eqref{p_k} is defined as
        \begin{align*}
            L(x, \phi_1, \dots, \phi_r) := f(x) - \sum_{i=1}^{r} g_i(x) \phi_i(x), \quad \phi_i(x) \in \Sigma^{l_i}[x]_{2(N-\tilde{w}_i)}.
        \end{align*}
        Here, the $l_i$-block circulant SOS polynomials $\phi_i$ is used.  
        Then, the Lagrangian dual problem for \eqref{p_k} can be expressed as
        \begin{align*}
            \max_{\phi_i}~\min_x ~ L(x, \phi_1, \dots, \phi_r).
        \end{align*}

        Now, we consider the following problem  with fixed $\phi_i(x)$:
         \begin{align}  \label{eq:FixedPhi0}
            \min_x ~ L(x, \phi_1, \dots, \phi_r),
        \end{align}
        which can be regarded as an unconstrained problem and  also be written as
        \begin{align} \label{eq:FixedPhi}
            \max  \ \gamma
            \mbox{ s.t. }  L(x, \phi_1, \dots, \phi_r) - \gamma \geq 0.
        \end{align}
        Then, we can apply $l_0$-block circulant SOS relaxation  \cite{zheng2022unconstrained} to \eqref{eq:FixedPhi0}, or equivalently to
        \eqref{eq:FixedPhi}:
        \begin{align} \label{eq:bcSOS}
            \max_{\gamma, \phi_i} ~ \gamma \quad \mathrm{s.t.} ~ L(x, \phi_1, \dots, \phi_r) - \gamma \in \Sigma^{l_0}[x]_{2N}.
        \end{align}
        In the problem \eqref{eq:bcSOS}, the constraint $L(x, \phi_1, \dots, \phi_r) - \gamma \in \Sigma^{l_0}[x]_{2N}$ is equivalent to the existence of some $\phi_0(x) \in \Sigma^{l_0}[x]_{2N}$ such that $L(x, \phi_1, \dots, \phi_r) - \gamma = \phi_0(x)$. Thus, 
        \eqref{p_k} can be relaxed to the following problem:
        \begin{maxi}
            {}  
            {\gamma} 
            {\tag{$\tilde{\mathbb{P}}_K$} \label{p_k_tilde}}  
            {}  
            \addConstraint{f(x) - \gamma}{\in \Gamma}{:= \left \{ \phi_0 + \sum_{i=1}^{r} g_i(x) \phi_i(x) \; \middle| \; \phi_0 \in \Sigma^{l_0}[x]_{2N}, \, \phi_i \in \Sigma^{l_i}[x]_{2(N-\tilde{w}_i)} \right \}.} 
        \end{maxi}

        By Assumption \ref{Assumption2}, \eqref{p_k_tilde} has a feasible solution. Moreover, the optimal value of \eqref{p_k_tilde} is not greater
         than the optimal value of \eqref{p_k}, based on the inclusion relations between the set of  block circulant SOS polynomials and the set of
          nonnegative polynomials \eqref{eq:Inclusion}, 
          i.e., $\max (\tilde{\mathbb{P}}_K) \le \min (\mathbb{P}_K)$.

          Next, we show that the relaxation problem \eqref{p_k_tilde} of \eqref{p_k} can be transformed into an equivalent problem using the T-SDP discussed in Section \ref{t-sdp}.
  With    a monomial vector  $[x]_N,$ 
        we define monomial tensors by
        \begin{gather*}
            \XC_0 := \fold_{l_0}([x]_N), \quad \XC_i := \fold_{l_i}([x]_{N-\tilde{w}_i}), \quad i=1, \dots, r.
        \end{gather*}
        We determine third-order tensors $\AC_\alpha \in \Real^{m_0 \times m_0 \times l_0}$, $\DC_{i \alpha} \in \Real^{m_i \times m_i \times l_i}$ such that
        \begin{align}
            \label{coefficient_tensor_obj}
            \XC_0 * \XC_0^\mathsf{T} &= \sum_{0 \le |\alpha| \le 2N} \AC_\alpha x^\alpha \\
            \label{coefficient_tensor_con}
            g_i(x) \XC_i * \XC_i^\mathsf{T} &= \sum_{0 \le |\alpha| \le 2N} \DC_{i \alpha} x^\alpha, \quad i=1, \dots, r.
        \end{align}
        If $l_0 = 1$ and $ l_i = 1$ in (\ref{mlsplitting}), then $\fold_1([x]_N) = [x]_N$ and $ \fold_1([x]_{N - \tilde{w}_i}) = [x]_{N - \tilde{w}_i}$ hold
        and the monomial tensors $\XC_0$ and $\XC_i$ become the monomial vectors, and the third-order tensors $\AC_\alpha, \DC_{i \alpha}$
        to be determined become
         the matrices $A_\alpha \in \Real^{m_0l_0 \times m_0l_0}$, $D_{i \alpha} \in \Real^{m_il_i \times m_il_i}$.
        By Theorem \ref{Theorem3}, the condition $f(x) - \gamma \in \Gamma$ in  \eqref{p_k_tilde} can be expressed as
        \begin{align} \label{eq:Relaxprob}
            f(x) - \gamma &= \left \langle \ZC_0, \XC_0 * \XC_0^\mathsf{T} \right \rangle + \sum_{i=1}^{r} g_i(x) \left \langle \ZC_i, \XC_i * \XC_i^\mathsf{T} \right \rangle  \nonumber \\
            & = \sum_{0 \le |\alpha| \le 2N} \langle \ZC_0, \AC_\alpha \rangle x^\alpha + \sum_{i=1}^{r} \sum_{0 \le |\alpha| \le 2N} \langle \ZC_i, \DC_{i \alpha} \rangle x^\alpha,
        \end{align}
        where $\ZC_0 \succeq_{\TC}\OC$, $\ZC_i \succeq_{\TC}\OC$ are semidefinite tensors. 
        Since the objective function $f(x)$ has no constant term ($f(0) = 0$), by
        comparing the coefficients  on both sides of \eqref{eq:Relaxprob} using \eqref{express-f}, we obtain
        \begin{align*}
            -\gamma &= \langle \ZC_0, \AC_0 \rangle + \sum_{i=1}^{r} \langle \ZC_i, \DC_{i 0} \rangle \\
            b_\alpha &= \langle \ZC_0, \AC_\alpha \rangle + \sum_{i=1}^{r} \langle \ZC_i, \DC_{i \alpha} \rangle, \quad 0 < |\alpha| \le 2N.
        \end{align*}
       Consequently, we obtain the T-SDP problem that is equivalent to problem \eqref{p_k_tilde}:
        \begin{maxi}
            {}  
            {-\langle \ZC_0, \AC_0 \rangle - \sum_{i=1}^{r} \langle \ZC_i, \DC_{i 0} \rangle} 
            {\tag{$\mathbb{Q}^{l_0, l_1, \dots, l_r}_N$} \label{qkn}}  
            {}  
            \addConstraint{\langle \ZC_0, \AC_\alpha \rangle + \sum_{i=1}^{r} \langle \ZC_i, \DC_{i \alpha} \rangle}{= b_\alpha, \ \quad 0 < |\alpha| \le 2N} 
            \addConstraint{\ZC_0 \succeq_{\TC}\OC, \  \ZC_i \succeq_{\TC}\OC, \quad i=1, \dots, r.} 
        \end{maxi}
        The above 
        \eqref{qkn} is an extension of the SDP relaxation problem derived from Parrilo's SOS relaxation
        \cite{parrilo2003semidefinite} to a T-SDP problem.
        The dual problem \ref{qknstar} of \ref{qkn} is expressed as
        \begin{mini}
            {}
            {\sum_{0 < |\alpha| \le 2N} b_\alpha y_\alpha} 
            {\tag*{($\mathbb{Q}^{l_0, l_1, \dots, l_r}_N)^*$} \label{qknstar}}  
            {}
            \addConstraint{\sum_{0 < |\alpha| \le 2N} \AC_\alpha y_\alpha}{\succeq_{\TC} -\AC_0} 
            \addConstraint{\sum_{0 < |\alpha| \le 2N} \DC_{i \alpha} y_\alpha}{\succeq_{\TC} -\DC_{i 0},}{\quad i=1, \dots, r.} 
        \end{mini}

        We notice that if we let $\MC_N^{l_0}(y) := \AC_0 + \sum_{0 < |\alpha| \le 2N} \AC_\alpha y_\alpha$, $\MC_{N-\tilde{w}_i}^{l_i}(g_i y) := \DC_{i 0} + \sum_{0 < |\alpha| \le 2N} \DC_{i \alpha} y_\alpha$, then the constraints of \ref{qknstar} can be written as $\MC_N^{l_0}(y) \succeq_{\TC}\OC$, $\MC_N^{l_i}(g_iy) \succeq_{\TC}\OC, \, i=1, \dots, r$.
        For the primal problem \eqref{qkn}, its dual \ref{qknstar} is an extension of the SDP relaxation problem derived by Lasserre's SDP relaxation \cite{lasserre2001global} to a T-SDP problem.
        The positive integers $l_0$ and $l_i$ determined to satisfy \eqref{mlsplitting} are the parameters in the block circulant SOS relaxation.
        By choosing a relaxation level $N$,  $l_0$ and $l_i$, the relaxation problems \eqref{qkn} and \ref{qknstar} are uniquely determined.

        From \eqref{mlsplitting}, (\ref{coefficient_tensor_obj}) and (\ref{coefficient_tensor_con}), the size of the T-SDP relaxation problem \eqref{qkn}
         is given with $\AC_\alpha, \ZC_0 \in \SymMat^{m_0 \times m_0 \times l_0}$ and $\DC_{i \alpha},\ZC_i \in \SymMat^{m_i \times m_i \times l_i}$. As discussed in Section \ref{t-sdp}, \eqref{qkn} can be transformed into an equivalent SDP of smaller size, depending on  even or odd $l_0, l_i$,
         which includes hermitian matrices of size determined by the sum of $\frac{l_0+2}{2}$ or $\frac{l_0+1}{2}$ of $m_0 \times m_0$ size
         and the sum of $\frac{l_i+2}{2}$ or $\frac{l_i+1}{2}$ of $m_i \times m_i$ size.
          In this case, we need to transform the SDP formulated with complex matrices to real-symmetric matrices \cite{aps2020mosek}.

        For the case of $l_0 = 1$ and $l_i = 1$, the T-SDP relaxation problem \eqref{qkn} coincides with the  SDP relaxation problem
        which involves  real symmetric matrices of size $m_0l_0$  
        and  $m_il_i$.  
        Therefore, the T-SDP relaxation problem can be reduced to a smaller problem than the  SDP relaxation in \cite{lasserre2001global}, as shown in Table~\ref{tab:size-and-number}.
        \begin{table}[htpb]
            \caption{Comparing the number and size of positive semidefinite (PSD) matrices for the basic
             SOS relaxation and the block circulant SOS relaxation
     }
            \label{tab:size-and-number}
            \begin{center}
                \scriptsize
                \begin{tabular}{cl} \hline
                    Relaxation & number and size of PSD matrices in the $N$th-level relaxation
                    \\ \hline
                   Basic SOS &
                    \begin{tabular}{l}
                        [ $1$*$(\binom{n + N}{N} \times \binom{n + N}{N})$ PSD matrix, \\
                        $1$*$(\sum_{i=1}^{r}$ $\binom{n + N - \tilde{w}_i}{N - \tilde{w}_i} \times \binom{n + N - \tilde{w}_i}{N - \tilde{w}_i})$ PSD matrices]
                    \end{tabular} \\
                    Block circulant SOS &
                    \begin{tabular}{l}
                        [($\frac{l_0 + 2}{2}$ or $\frac{l_0 + 1}{2}$) *$\left(\frac{\binom{n + N}{N}}{l_0} \times \frac{\binom{n + N}{N}}{l_0}\right)$ PSD matrices, \\
                        $\left(\sum_{i=1}^{r} \frac{l_i + 2}{2} \ \textrm{or} \ \sum_{i=1}^{r} \frac{l_i + 1}{2}\right)$* $\left(\frac{\binom{n + N - \tilde{w}_i}{N - \tilde{w}_i}}{l_i} \times \frac{\binom{n + N - \tilde{w}_i}{N - \tilde{w}_i}}{l_i}\right)$ PSD matrices]
                    \end{tabular} \\
                    \hline
                \end{tabular}
            \end{center}
        \end{table}

        \subsection{Feasibility and global optimality}
        \label{feas_and_opti}
        We present theoretical analysis of the T-SDP relaxation problems \eqref{qkn} and \ref{qknstar} for the constrained POPs proposed in Section~\ref{tsdp-relax}.
 	We also show that the sequence of relaxation problems $\left\{ (\mathbb{Q}_N^{l_0, l_1, \dots, l_r})^* \right\}$
	 for a relaxation level $N$ converges to the optimal value of \eqref{p_k} when $N \to \infty$.

       Using the relation between block circulant SOS polynomials and semidefinite block circulant matrices described in Section \ref{bcircSOSpoly}, we present the following theorem, which serves as  a necessary and sufficient condition for the feasibility  of \eqref{qkn}.

        \begin{theorem} \label{Theorem4}
         Let $f:\Real^n \to \Real$ be a real-valued polynomial of degree $d$ with the zero constant term. 
         Assume that there exist the positive integers $m_0, l_0, m_i, l_i,  i=1, \dots, r$ such that
          $s(N) = m_0 l_0$, $s(N-\tilde{w}_i) = m_i l_i$. Then \eqref{qkn} has a feasible solution if and only if there exists some $\gamma \in \Real$ and third-order tensors $\ZC_0 \in \SymMat^{m_0 \times m_0 \times l_0}_+$, $\ZC_i \in \SymMat^{m_i \times m_i \times l_i}_+,  i=1, \dots, r$ such that
        \[ f(x) - \gamma = \left \langle \ZC_0, \XC_0 * \XC_0^\mathsf{T} \right \rangle + \sum_{i=1}^{r} g_i(x) \left \langle \ZC_i, \XC_i * \XC_i^\mathsf{T} \right \rangle. \]
 	\end{theorem}
        \textit{Proof.} Assume that $(\ZC_0, \ZC_1, \dots, \ZC_r)$ is a feasible solution of \eqref{qkn}. More precisely,
        \begin{equation}
            \begin{aligned}
                \langle \ZC_0, \AC_\alpha \rangle + \sum_{i=1}^{r} \langle \ZC_i, \DC_{i \alpha} \rangle &= b_\alpha, \quad 0 < |\alpha| \le 2N \\
                \ZC_0 &\succeq_{\TC}\OC, \
                \ZC_i \succeq_{\TC}\OC, \quad i=1, \dots, r.
            \end{aligned} \notag
        \end{equation}
        Since the objective function in \eqref{p_k} has no constant term,  we have
        \begin{align*}
            f(x) &= \sum_{0 < |\alpha| \le 2N} b_\alpha x^\alpha
            = \sum_{0 < |\alpha| \le 2N} \left( \langle \ZC_0, \AC_\alpha \rangle x^\alpha + \sum_{i=1}^{r} \langle \ZC_i, \DC_{i \alpha} \rangle x^\alpha \right) \\
            &= \sum_{0 < |\alpha| \le 2N} \langle \ZC_0, \AC_\alpha \rangle x^\alpha + \sum_{i=1}^{r} \sum_{0 < |\alpha| \le 2N} \langle \ZC_i, \DC_{i \alpha} \rangle x^\alpha \\
            &= \left \langle \ZC_0, \XC_0 * \XC_0^\mathsf{T} \right \rangle - \langle \ZC_0, \AC_0 \rangle + \sum_{i=1}^{r} \left( g_i(x) \left \langle \ZC_i, \XC_i * \XC_i^\mathsf{T} \right \rangle - \langle \ZC_i, \DC_{i 0} \rangle \right),
        \end{align*}
        where  \eqref{coefficient_tensor_obj} and \eqref{coefficient_tensor_con} are used in the last equality. 
        Since $\AC_0 \in \Real^{m_0 \times m_0 \times l_0}$, $\DC_{i 0} \in \Real^{m_i \times m_i \times l_i}$ are all 0 except for the $(1, 1, 1)$th
         element, which is 1 by its definition,
        \[ f(x) + [\ZC_0]_{(1,1,1)} + \sum_{i=1}^{r} [\ZC_i]_{(1,1,1)} = \left \langle \ZC_0, \XC_0 * \XC_0^\mathsf{T} \right \rangle + \sum_{i=1}^{r} g_i(x) \left \langle \ZC_i, \XC_i * \XC_i^\mathsf{T} \right \rangle. \]
        The desired result follows by setting $\gamma = - [\ZC_0]_{(1,1,1)} - \sum_{i=1}^{r} [\ZC_i]_{(1,1,1)}$.

        On the other hand, assume that $\gamma \in \Real$ and $\ZC_0 \in \SymMat^{m_0 \times m_0 \times l_0}_+$, $\ZC_i \in \SymMat^{m_i \times m_i \times l_i}_+, \, i=1, \dots, r$ exist such that 
        \[ f(x) - \gamma = \left \langle \ZC_0, \XC_0 * \XC_0^\mathsf{T} \right \rangle + \sum_{i=1}^{r} g_i(x) \left \langle \ZC_i, \XC_i * \XC_i^\mathsf{T} \right \rangle. \]
        Then, from (\ref{coefficient_tensor_obj}) and (\ref{coefficient_tensor_con}), we have
        \begin{align*}
            f(x) = \sum_{0 < |\alpha| \le 2N} \langle \ZC_0, \AC_\alpha \rangle x^\alpha + \langle \ZC_0, \AC_0 \rangle + \sum_{0 < |\alpha| \le 2N} \sum_{i=1}^{r} \langle \ZC_i, \DC_{i \alpha} \rangle x^\alpha + \sum_{i=1}^{r} \langle \ZC_i, \DC_{i 0} \rangle + \gamma.
        \end{align*}
        Therefore, by comparing the coefficients for each monomial on both sides, we obtain
        \begin{align*}
            \langle \ZC_0, \AC_0 \rangle + \sum_{i=1}^{r} \langle \ZC_i, \DC_{i 0} \rangle + \gamma &= 0 \\
            \langle \ZC_0, \AC_\alpha \rangle + \sum_{i=1}^{r} \langle \ZC_i, \DC_{i \alpha} \rangle &= b_\alpha, \quad 0 < \alpha \le 2N.
        \end{align*}
        Thus, $\ZC_0 \in \SymMat^{m_0 \times m_0 \times l_0}_+$, $\ZC_i \in \SymMat^{m_i \times m_i \times l_i}_+, \, i=1, \dots, r$ are feasible solutions of \eqref{qkn}. $\square$

        We show the relation between the SDP relaxation problem by Parrilo \cite{parrilo2003semidefinite} and the T-SDP relaxation problem of the proposed method in the following theorem. In particular, we discuss the necessary and sufficient conditions under which the optimal values of SDP relaxation and T-SDP relaxation are equivalent.

        \begin{theorem} \label{Theorem5}
        Let $(\mathbb{Q}_N^{1,1,\dots,1})$ be the SDP relaxation problem for \eqref{p_k} with  $l_0 = 1$, $l_1 = 1, \dots, l_r = 1$
         and $(\mathbb{Q}_N^{l_0,l_1,\dots,l_r})$ be the T-SDP relaxation problem with  $l_0$, $l_1, \dots, l_r$ such that $\max\{l_0, l_1, \dots, l_r\} \ge 2$. Let the optimal values of the relaxation problem be $\max (\mathbb{Q}_N^{1,1,\dots,1})$, $\max (\mathbb{Q}_N^{l_0,l_1,\dots,l_r})$, respectively. Then $\max (\mathbb{Q}_N^{1,1,\dots,1}) = \max (\mathbb{Q}_N^{l_0,l_1,\dots,l_r})$ if and only if there exist optimal solutions $Z_0^*$, $Z_i^*$ of $(\mathbb{Q}_N^{1,1,\dots,1})$ 
        represented by $l_0$-block circulant and $l_i$-block circulant matrices, respectively.
        \end{theorem}

        \textit{Proof.} From Theorem~\ref{Theorem4}, 
        the constraints of $(\mathbb{Q}_N^{1,1,\dots,1})$ is equivalent to
        \[ f(x) - \gamma = \left\langle Z_0, [x]_N * [x]_N^\mathsf{T} \right\rangle + \sum_{i=1}^{r} g_i(x)\left\langle Z_i, [x]_{N-\tilde{w}_i} * [x]_{N-\tilde{w}_i}^\mathsf{T} \right\rangle, \]
        and the constraints of $(\mathbb{Q}_N^{l_0,l_1,\dots,l_r})$ is equivalent to
        \[ f(x) - \gamma = \left\langle \ZC_0, \XC_0 * \XC_0^\mathsf{T} \right\rangle + \sum_{i=1}^{r} g_i(x)\left\langle \ZC_i, \XC_i * \XC_i^\mathsf{T} \right\rangle. \]
        Then,
        \begin{align*}
            &\quad\left\langle \ZC_0, \XC_0 * \XC_0^\mathsf{T} \right\rangle + \sum_{i=1}^{r} g_i(x)\left\langle \ZC_i, \XC_i * \XC_i^\mathsf{T} \right\rangle \\
            &= \frac{1}{l_0}\mathrm{Tr}\left( \bcirc(\XC_0)^\mathsf{T} \bcirc(\ZC_0) \bcirc(\XC_0) \right) + \sum_{i=1}^{r} g_i(x) \frac{1}{l_i}\mathrm{Tr}\left( \bcirc(\XC_i)^\mathsf{T} \bcirc(\ZC_i) \bcirc(\XC_i) \right) \\
                        &= \frac{1}{l_0} \left\langle \bcirc(\XC_0), \bcirc(\ZC_0) \bcirc(\XC_0) \right\rangle + \sum_{i=1}^{r} g_i(x) \frac{1}{l_i} \left\langle \bcirc(\XC_i), \bcirc(\ZC_i) \bcirc(\XC_i) \right\rangle
                        \end{align*}
            \begin{align*}
       &  \hspace{-8mm} =    \left\langle \unfold(\XC_0), \bcirc(\ZC_0) \unfold(\XC_0) \right\rangle + \sum_{i=1}^{r} g_i(x) \left\langle \unfold(\XC_i), \bcirc(\ZC_i) \unfold(\XC_i) \right\rangle \\
            & \hspace{-8mm} = \left\langle \bcirc(\ZC_0), [x]_N * [x]_N^\mathsf{T} \right\rangle + \sum_{i=1}^{r} g_i(x) \left\langle \bcirc(\ZC_i), [x]_{N-\tilde{w}_i} * [x]_{N-\tilde{w}_i}^\mathsf{T} \right\rangle.
        \end{align*}
        Hence, we have $\max (\mathbb{Q}_N^{1,1,\dots,1}) \ge \max (\mathbb{Q}_N^{l_0,l_1,\dots,l_r})$ since a feasible solution $(\ZC_0^*, \ZC_1^*, \dots, \ZC_r^*)$ of $(\mathbb{Q}_N^{l_0,l_1,\dots,l_r})$ becomes a feasible solution of $(\mathbb{Q}_N^{1,1,\dots,1})$.

        Now, we assume that $\max (\mathbb{Q}_N^{1,1,\dots,1}) = \max (\mathbb{Q}_N^{l_0,l_1,\dots,l_r})$. Let $(\ZC_0^*, \ZC_1^*, \dots, \ZC_r^*)$ be an optimal solution of $(\mathbb{Q}_N^{l_0,l_1,\dots,l_r})$.
        Then $(\bcirc(\ZC_0^*), \bcirc(\ZC_1^*), \dots, \bcirc(\ZC_r^*))$ is an optimal solution of $(\mathbb{Q}_N^{1,1,\dots,1})$. Clearly, $\bcirc(\ZC_0^*)$ and $\bcirc(\ZC_i^*)$ are $l_0$-block circulant and $l_i$-block circulant matrices, respectively. On the other hand, we can take $\bcirc_{l_0}^{-1}(Z_0^*), \bcirc_{l_i}^{-1}(Z_i^*)$, assuming that there exists an optimal solution $Z_0^*$ and $Z_i^*$ of $(\mathbb{Q}_N^{1,1,\dots,1})$
         represented as $l_0$-block circulant and $l_i$-block circulant matrices, respectively, which is also
          a feasible solution of $(\mathbb{Q}_N^{l_0,l_1,\dots,l_r})$.
        Consequently, $\max (\mathbb{Q}_N^{1,1,\dots,1}) \le \max (\mathbb{Q}_N^{l_0,l_1,\dots,l_r})$ holds, and together with $\max (\mathbb{Q}_N^{1,1,\dots,1}) \ge \max (\mathbb{Q}_N^{l_0,l_1,\dots,l_r})$, we have $\max (\mathbb{Q}_N^{1,1,\dots,1}) = \max (\mathbb{Q}_N^{l_0,l_1,\dots,l_r})$. $\square$

        Now, we discuss the global optimality of the T-SDP relaxation problem and  the convergence of the relaxation problem sequence $\left\{ \text{\ref{qknstar}} \right\}$ with respect to the relaxation level $N$, extending the Lasserre's hierarchy in \cite{lasserre2001global}.
 
        \begin{theorem} \label{Theorem6}
        Let $f:\Real^n \to \Real$ be a real-valued polynomial of degree $d$ and $K$ be the compact set. Let Assumption~\ref{Assumption2}
         hold, and let $p_K^* := \min_{x \in K} f(x)$.
         Let $x^*$ be a global optimal solution of \eqref{p_k} and
         \[ y^* = [x_1^*, \dots, x_n^*, (x_1^*)^2, x_1^*x_2^*, \dots, (x_n^*)^2, \dots, (x_1^*)^{2N}, \dots, (x_n^*)^{2N} ]^\mathsf{T}. \]
         Then, \\
\noindent
        (a) For fixed parameters $l_0, l_i$,  we have 
        \[ \inf (\mathbb{Q}_N^{l_0,l_1,\dots,l_r}) \, \uparrow \, p_K^*, \]
         as $N \to \infty$.
        Moreover, for $N$ sufficiently large, there is no duality gap between \eqref{qkn} and its dual {\upshape \ref{qknstar}} if $K$ has a nonempty interior. \\
\noindent
        (b) If $f(x) - p_K^*$ can be represented in the form (\ref{assumption2}), i.e.,
        \begin{align*}
            f(x) - p_K^* = q_0(x) + \sum_{i=1}^{r} g_i(x) q_i(x)
        \end{align*}
        for a block circulant SOS polynomial $q_0(x)$ of degree at most $2N$, and some  block circulant SOS polynomials $q_i(x)$ of degree at most 
        $2(N-\tilde{w}_i)$, 
        then
        \[ \min \text{\upshape \ref{qknstar}} = p_K^* = \max \text{\upshape (\ref{qkn})} \] and 
        $y^*$
        is a global minimizer of {\upshape \ref{qknstar}}.
        \end{theorem}
        \textit{Proof.} (a) Let $x^*$ be a global optimal solution of \eqref{p_k}, and
        \[ y^* = [x_1^*, \dots, x_n^*, (x_1^*)^2, x_1^*x_2^*, \dots, (x_n^*)^2, \dots, (x_1^*)^{2N}, \dots, (x_n^*)^{2N} ]^\mathsf{T}. \]
        Then, $\MC_N^{l_0}(y^*) = \fold_{l_0}(y^*)*\left( \fold_{l_0}(y^*) \right)^\mathsf{T} $, $\MC_{N-\tilde{w}_i}^{l_i}(g_iy^*) = g_i(x^*) \fold_{l_i}\left(y^*\right)*\fold_{l_i}\left(y^*\right)^\mathsf{T}, \, i=1, \dots, r$ are semidefinite tensors.
        Moreover, since $b_\alpha$ are the coefficients corresponding to the monomials $x^\alpha$, $\sum_{0 < \alpha \le 2N} b_\alpha y_\alpha^*$ is
        equal to $p_K^*$. Therefore, $y^*$ is a feasible solution of \ref{qknstar} with the
        objective function value $p_K^*$, hence, $\inf \textrm{\upshape \ref{qknstar}} \le p_K^*$.

        Now, for fixed $l_0, l_1, \dots, l_r$, we consider any $N' \ge N$ such that $s(N) = m_0l_0, s(N-\tilde{w}_i) = m_il_i$ and $s(N') = m_0'l_0, s(N'-\tilde{w}_i) = m_i'l_i$.
        In this case, since $s(N') - s(N)$ and $s(N'-\tilde{w}_i) - s(N-\tilde{w}_i)$ obviously have factors $l_0$ and $l_i$, respectively, it is possible to
        represent $\MC_{N}^{l_0}(y)$ and $\MC_{N-\tilde{w}_i}^{l_i}(g_iy)$ to be subtensors of $\MC_{N'}^{l_0}(y)$ and $\MC_{N'-\tilde{w}_i}^{l_i}(g_iy)$, respectively, by arranging the monomial vectors appropriately when creating the T-SDP relaxation problem.
        Specifically, if $\MC_{N'}^{l_0}(y) \succeq_{\TC}\OC$ and $\MC_{N'-\tilde{w}_i}^{l_i}(g_iy) \succeq_{\TC}\OC$,
        then  $\MC_{N}^{l_0}(y)$, $\MC_{N-\tilde{w}_i}^{l_i}(g_iy)$ can be formed
        such that $\MC_{N}^{l_0}(y) \succeq_{\TC}\OC$, $\MC_{N-\tilde{w}_i}^{l_i}(g_iy) \succeq_{\TC}\OC$, respectively (the details are described in the subsequent discussion and Example~\ref{Example2}).
        Thus, for any solution $y$ of $(\mathbb{Q}_{N'}^{l_0,l_1,\dots,l_r})^*$, the adequate truncated vector $y'$ is a feasible solution of \ref{qknstar}.
        Since $2N' \ge 2N \ge d$, $b_\alpha$ corresponding to $|\alpha| > \binom{n + d}{d}$ is $0$, therefore, $\sum_{0 < |\alpha| \le 2N} b_\alpha y_\alpha = \sum_{0 < |\alpha| \le 2N'} b_\alpha y_\alpha'$ holds. Then the objective value in the feasible solution $y$ of \ref{qknstar} is equal to 
        that in the feasible solution $y'$ of $(\mathbb{Q}_{N'}^{l_0,l_1,\dots,l_r})^*$. Therefore, $\inf (\mathbb{Q}_{N'}^{l_0,l_1,\dots,l_r})^* \ge \inf \text{\ref{qknstar}}$ for $N' \ge N$.

        From Assumption \ref{Assumption2},
        \[ f(x) - p_K^* + \epsilon = q_0(x) + \sum_{i=1}^{r} g_i(x) q_i(x), \]
        where $q_0(x)$ is an $l_0$-block circulant SOS polynomial of degree $2N$ and $q_i(x)$ is  an $l_i$-block circulant SOS polynomials
        of degree $2(N-\tilde{w}_i)$. In addition, from the definition of the block circulant SOS polynomial, there exist $t_0$ and $t_k, \, k=1, \dots, r$,
        $q_0$, and $q_i$ that can be expressed as
        \begin{gather*}
            q_0(x) = \sum_{j=1}^{t_0} \sum_{k=1}^{l_0} \left((u_{jk}^0)^\mathsf{T}[x]_N \right)^2 \\
            q_i(x) = \sum_{j=1}^{t_i} \sum_{k=1}^{l_i} \left((u_{jk}^i)^\mathsf{T}[x]_{N-\tilde{w}_i}\right)^2, \quad i=1, \dots, r.
        \end{gather*}
        We define
        \begin{align*}
            \ZC_0 := \bcirc_{l_0}^{-1}(\sum_{j=1}^{t_0} \sum_{k=1}^{l_0} u_{jk}^0(u_{jk}^0)^\mathsf{T}), \quad \ZC_i := \bcirc_{l_i}^{-1}(\sum_{j=1}^{t_i} \sum_{k=1}^{l_i} u_{jk}^i(u_{jk}^i)^\mathsf{T}), \quad i=1, \dots, r.
        \end{align*}
        Then, from Theorem \ref{Theorem1}, we have $\ZC_0, \ZC_k \succeq_{\TC}\OC$. Since
        \begin{align*}
            f(x) - p^*_K + \epsilon &= q_0(x) + \sum_{i=1}^{r} g_i(x)q_i(x) \\
            &= \left\langle \sum_{j=1}^{t_0} \sum_{k=1}^{l_0} u_{jk}^0(u_{jk}^0)^\mathsf{T}, [x]_N [x]_N^\mathsf{T} \right\rangle + \sum_{i=1}^{r} g_i(x)\left\langle \sum_{j=1}^{t_i} \sum_{k=1}^{l_i} u_{jk}^i(u_{jk}^i)^\mathsf{T}, [x]_{N-\tilde{w}_i} [x]_{N-\tilde{w}_i}^\mathsf{T} \right\rangle \\
            &= \left\langle \ZC_0, \XC_0 * \XC_0^\mathsf{T} \right\rangle + \sum_{i=1}^{r} g_i(x)\left\langle \ZC_i, \XC_i * \XC_{i}^\mathsf{T} \right\rangle,
        \end{align*}
        holds from Theorem \ref{Theorem4}, $(\ZC_0, \ZC_1, \dots, \ZC_r)$ is a feasible solution of \eqref{qkn} and $[\ZC_0]_{(1,1,1)} + \sum_k [\ZC_k]_{(1,1,1)} = -(p_K^* - \epsilon)$. Thus, we obtain
        \begin{align}
        p_K^* - \epsilon \le \sup (\text{\ref{qkn}}) \le \inf \text{\ref{qknstar}} \le p_K^*.
            \label{eq:boundness-of-qknstar}
        \end{align}

        Next, we prove that 
        there is no duality gap between \eqref{qkn} and its dual \ref{qknstar} for any $N$ such that $ N \ge N_0$ where $N_0$ is an initial relaxation level.
        Let $\mu$ be a probability measure with uniform distribution in $K$, which has a strictly positive density $f$ with respect to the Lebesgue
        measure and satisfies that
        \[ y_\alpha = \int x^\alpha d\mu < +\infty, \quad \text{for all} \;  \alpha \ \text{such that} \ 0 \le |\alpha| \le 2N_0. \]
        Then, from the discussion in Theorem 6 in \cite{zheng2022unconstrained}, it follows that $\MC_N^{l_0}(y_\mu) \succ_{\TC}\OC$ and
        $\MC_{N-\tilde{w}_i}^{l_i}(g_iy_\mu) \succ_{\TC}\OC, \, i=1, \dots, r$ for $y_\mu = \{ y_\alpha \}$, which is a feasible interior point of \ref{qknstar}.
        In addition, \ref{qknstar} is bounded below by \eqref{eq:boundness-of-qknstar}.
        By strong duality  of T-SDP in Theorem~\ref{Theorem2},
        this indicates that \eqref{qkn} is solvable, together with \eqref{qkn} is feasible,
         and there is
        no duality gap between \eqref{qkn} and \ref{qknstar}.

        (b) If $f(x) - p^*_K$ can be expressed as in (\ref{assumption2}), then  we can construct matrices
        $\ZC_0 \succeq_{\TC}\OC $, $\ZC_i \succeq_{\TC}\OC, \, i=1, \dots, r$
        to be a feasible solution of \eqref{qkn} whose objective value is $-[\ZC_0]_{(1,1,1)} - \sum_k [\ZC_k]_{(1,1,1)} = p^*_K$
        with polynomials $q_0(x)$ and $q_i(x)$ of degree at most $2N$ and $2(N-\tilde{w}_i)$, respectively,
          as in the proof of (a).
         As a result,
         from $p_K^* \le \sup (\text{\ref{qkn}}) \le \inf \text{\ref{qknstar}} \le p_K^* $, we have $\max (\text{\ref{qkn}}) = p^*_K = \min \text{\ref{qknstar}}$,
         and $(\ZC_0, \ZC_1, \dots, \ZC_k)$ is an optimal solution of \eqref{qkn}. Furthermore, since $b_\alpha$ are the coefficients corresponding to each
         monomial, we clearly see that $y^*$ is a global optimal solution of \ref{qknstar}. $\square$

        For any two third-order tensors $\AC \in \Real^{m \times n \times l}$ and $\BC \in \Real^{m' \times n' \times l}$, $\AC$ is called
        a subtensor of $\BC$ if each frontal slice $\AC^{(i)}$ of $\AC$ is a principal submatrix of each frontal slice $\BC^{(i)}$ of $\BC$, respectively.

        We illustrate the proof of (a) in Theorem~\ref{Theorem6} with the following examples.
        \begin{example} \label{Example2}
        We let $n=2$, the relaxation level $N = 2$ and $N' = 3$. We fix $l_0 = 2$.
         \end{example}
        In this case, $s(N) = 6$. We derive T-SDP relaxation by a monomial vector
        \begin{align*}
            [x]_2 = [1, x_1, x_2, x_1^2, x_1x_2, x_2^2]^\mathsf{T}.
        \end{align*}
        Then $\MC_2^2(y)$ is a third-order tensor of size $3 \times 3 \times 2$ whose frontal slices are given by
        \begin{align*}
           \MC_2^2(y)^{(1)} = \left[ \begin{array}{@{\,}ccc@{\,}}
                                           1 + y_{40} & y_{10} + y_{31} & y_{01} + y_{22} \\
                                           y_{10} + y_{31} & y_{20} + y_{22} & y_{11} + y_{13} \\
                                           y_{01} + y_{22} & y_{11} + y_{13} & y_{02} + y_{04}
                                       \end{array} \right], \quad
           \MC_2^2(y)^{(2)} = \left[ \begin{array}{@{\,}ccc@{\,}}
                                           2y_{20} & y_{11} + y_{30} & y_{02} + y_{21} \\
                                           y_{11} + y_{30} & 2y_{21} & 2y_{12} \\
                                           y_{21} + y_{02} & 2y_{12} & 2y_{03}
                                       \end{array} \right].
        \end{align*}
        With $s(N') = 10$,  T-SDP relaxation can be  derived using a monomial vector
        \begin{align*}
            [x]_3 = [1, x_1, x_2, x_1^2, x_1x_2, x_2^2, x_1^3, x_1^2x_2, x_1x_2^2, x_2^3]^\mathsf{T}.
        \end{align*}
        Then, $\MC_3^2(y)$ is a third-order tensor of size $5 \times 5 \times 2$ whose frontal slices are given by
        \begin{align*}
           \MC_3^2(y)^{(1)} &= \left[ \begin{array}{@{\,}ccccc@{\,}}
                                           1 + y_{04} & y_{10} + y_{32} & y_{01} + y_{23} & y_{20} + y_{14} & y_{11} + y_{05} \\
                                           y_{10} + y_{32} & y_{20} + y_{60} & y_{11} + y_{51} & y_{30} + y_{42} & y_{21} + y_{33} \\
                                           y_{01} + y_{23} & y_{11} + y_{51} & y_{02} + y_{42} & y_{21} + y_{33} & y_{12} + y_{24} \\
                                           y_{20} + y_{14} & y_{30} + y_{42} & y_{21} + y_{33} & y_{40} + y_{24} & y_{32} + y_{15} \\
                                           y_{11} + y_{05} & y_{21} + y_{33} & y_{12} + y_{24} & y_{31} + y_{15} & y_{22} + y_{06}
                                       \end{array} \right], \\
            \\
           \MC_3^2(y)^{(2)} &= \left[ \begin{array}{@{\,}ccccc@{\,}}
                                           2y_{02} & y_{12} + y_{30} & y_{03} + y_{21} & y_{12} + y_{22} & y_{30} + y_{13} \\
                                           y_{12} + y_{30} & 2y_{40} & 2y_{31} & y_{22} + y_{50}  & y_{13} + y_{42} \\
                                           y_{03} + y_{21} & 2y_{31} & 2y_{22} & y_{13} + y_{41} & y_{04} + y_{32} \\
                                           y_{12} + y_{22} & y_{22} + y_{50} & y_{13} + y_{42} & 2y_{32}  & 2y_{23} \\
                                           y_{03} + y_{13} & y_{13} + y_{42} & y_{04} + y_{32} & 2y_{23} & 2y_{14}
                                       \end{array} \right].
        \end{align*}
        Obviously, $\MC_2^2(y)$ is not a subtensor of $\MC_3^2(y)$. However, if we derive the T-SDP relaxation problem by a monomial vector
        \begin{align*}
            [x]_3 = [1, x_1, x_2, x_1^3, x_1^2x_2, x_1^2, x_1x_2, x_2^2, x_1x_2^2, x_2^3]^\mathsf{T},
        \end{align*}
        then each frontal slice of $\MC_3^2(y)$ is given by
        \begin{align*}
           \MC_3^2(y)^{(1)} &= \left[ \begin{array}{@{\,}ccc:cc@{\,}}
                                           1 + y_{40} & y_{10} + y_{31} & y_{01} + y_{22} & y_{30} + y_{32} & y_{21} + y_{23} \\
                                           y_{10} + y_{32} & y_{20} + y_{22} & y_{11} + y_{13} & y_{40} + y_{23} & y_{31} + y_{14} \\
                                           y_{01} + y_{22} & y_{11} + y_{13} & y_{02} + y_{04} & y_{31} + y_{14} & y_{22} + y_{05} \\ \hdashline
                                           y_{30} + y_{32} & y_{40} + y_{23} & y_{31} + y_{14} & y_{24} + y_{60} & y_{51} + y_{15} \\
                                           y_{21} + y_{23} & y_{31} + y_{14} & y_{22} + y_{05} & y_{51} + y_{15} & y_{42} + y_{06}
                                       \end{array} \right], \\
            \\
           \MC_3^2(y)^{(2)} &= \left[ \begin{array}{@{\,}ccc:cc@{\,}}
                                           2y_{20} & y_{11} + y_{30} & y_{02} + y_{21} & y_{12} + y_{50} & y_{03} + y_{41} \\
                                           y_{11} + y_{30} & 2y_{21} & 2y_{12} & y_{22} + y_{41}  & y_{13} + y_{32} \\
                                           y_{02} + y_{21} & 2y_{12} & 2y_{03} & y_{13} + y_{432} & y_{04} + y_{23} \\ \hdashline
                                           y_{12} + y_{50} & y_{22} + y_{41} & y_{13} + y_{32} & 2y_{42}  & 2y_{33} \\
                                           y_{03} + y_{41} & y_{13} + y_{32} & y_{04} + y_{23} & 2y_{33} & 2y_{24}
                                       \end{array} \right],
        \end{align*}
        where $\MC_2^2(y)$ is a subtensor of $\MC_3^2(y)$.

 	Example \ref{Example2} illustrates that
        for any $N$ and $N'$ such that $N' \ge N$, if $l_0$ is fixed, $\MC_N^{l_0}(y)$ can be formed to be a subtensor of $\MC_{N'}^{l_0}(y)$ by arranging the monomial vector $[x]_{N'}$ for the derivation of T-SDP relaxation with relaxation level
         $N'$ such that the monomial tensor $\fold_{l_0}([x]_N)$ is a subtensor of $\fold_{l_0}([x]_{N'})$.
        Similarly, if $l_i$ corresponding to each constraint is fixed, $\MC_N^{l_i}(g_iy)$ can be formed to be a subtensor of $\MC_{N'}^{l_i}(g_iy)$ by appropriately arranging the monomial vectors $[x]_{N'-\tilde{w}_i}$.
        For example, consider a constrained POP in two variables ($n=2$) with  the objective function of degree $d = 6$, the number of constraints $r = 2$, and degree $w_i = 1, \, i=1, 2$ for each constraint.
        Let $(\mathbb{Q}_N^{5,2,2})^*$ denote the T-SDP relaxation obtained with the relaxation level $N = 3$.
        Then $s(N) = 10$, $s(N-\tilde{w}_i) = 6$, $l_0 = 5$ and $l_i = 2$.
        Furthermore, if the relaxation level $N' = 4$, then $s(N') = 15$, $s(N' - \tilde{w}_i) = 10$, and if $m_0' = 3$, $m_i' = 5$,
        then $s(N') = m_0' l_0$ and $s(N' - \tilde{w}_i) = m_i' l_i$.
        Thus, we can also derive the T-SDP  relaxation with $l_0 = 5$ and $l_i = 2$,  denoted by $(\mathbb{Q}_{N'}^{5,2,2})^*$. In (a) of Theorem
        \ref{Theorem6}, we have claimed that the relation $\inf (\mathbb{Q}_{N}^{5,2,2})^* \le \inf (\mathbb{Q}_{N'}^{5,2,2})^*$ holds.
        Conversely, if $N'' = 5$, then $s(N'') = 21$ and $s(N'' - \tilde{w}_i) = 15$, but there are no positive integers $m_0''$ and $m_i''$ such that $s(N'') = m_0''l_0$ and $s(N'' - \tilde{w}_i) = m_i''l_i$.
        Therefore, any T-SDP relaxation with relaxation level $N'' = 5$ has no relation to $(\mathbb{Q}_{N}^{5,2,2})^*$ or $(\mathbb{Q}_{N'}^{5,2,2})^*$.
        Clearly, if we fix $l_0 = 1$ and $l_i = 1$, then $\inf (\mathbb{Q}_{N}^{1,1,1})^* \le \inf (\mathbb{Q}_{N'}^{1,1,1})^*$ for any $N' \ge N$.

        We have mentioned that   ``if $\MC_{N'}^{l_0}(y) \succeq_{\TC}\OC$ and $\MC_{N'-\tilde{w}_i}^{l_i}(g_iy) \succeq_{\TC}\OC$, then  $\MC_{N}^{l_0}(y)$, $\MC_{N-\tilde{w}_i}^{l_i}(g_iy)$ can be constructed
         such that $\MC_{N}^{l_0}(y) \succeq_{\TC}\OC$, $\MC_{N-\tilde{w}_i}^{l_i}(g_iy) \succeq_{\TC}\OC$, respectively"  in the proof of (a) of Theorem \ref{Theorem6}.
         As discussed earlier, $\MC_{N}^{l_0}(y)$, $\MC_{N-\tilde{w}_i}^{l_i}(g_iy)$ can be constructed
          to be subtensors of $\MC_{N'}^{l_0}(y)$, $\MC_{N'-\tilde{w}_i}^{l_i}(g_iy)$, respectively, in the T-SDP relaxation with  relaxation
           level $N$ such that $N' \ge N$, by  appropriately arranging the monomial vectors. Therefore, from the definition of subtensor, $\bcirc\left(\MC_{N}^{l_0}(y)\right)$ and $\bcirc\left(\MC_{N-\tilde{w}_i}^{l_i}(g_iy)\right)$ are
         the principal submatrices of $\bcirc\left(\MC_{N'}^{l_0}(y)\right)$ and $\bcirc\left(\MC_{N'-\tilde{w}_i}^{l_i}(g_iy)\right)$, respectively.
        Now, if we assume $\MC_{N'}^{l_0}(y) \succeq_{\TC}\OC$ and $\MC_{N'-\tilde{w}_i}^{l_i}(g_iy) \succeq_{\TC}\OC$, then
             we obtain $\bcirc\left(\MC_{N'}^{l_0}(y)\right) \succeq O$, $\bcirc\left(\MC_{N'-\tilde{w}_i}^{l_i}(g_iy)\right) \succeq O$ 
                  from Theorem \ref{Theorem1}. 
        A matrix $A$ is positive semidefinite if and only if the determinants of all principal minors of $A$ are nonnegative, so that
        the determinant of all principal minors of $\bcirc\left(\MC_{N}^{l_0}(y)\right)$, $\bcirc\left(\MC_{N-\tilde{w}_i}^{l_i}(g_iy)\right)$ also are nonnegative, respectively. Therefore, $\bcirc\left(\MC_{N}^{l_0}(y)\right) \succeq O$ and $\bcirc\left(\MC_{N-\tilde{w}_i}^{l_i}(g_iy)\right) \succeq O$. Thus, by Theorem \ref{Theorem1},  $\MC_{N}^{l_0}(y) \succeq\OC$ and $\MC_{N-\tilde{w}_i}^{l_i}(g_iy) \succeq\OC$ follows.

%
     \section{Numerical experiments}
        \label{expe}
        We compare 
         the basic SOS relaxation \cite{kim2005generalized} with the  T-SDP relaxation
        proposed in Section \ref{tsdp-relax} for constrained POPs  with ten test problems and show
        that the proposed T-SDP relaxation is more efficient than the basic SOS relaxation.
        The test problems are presented in detail in Appendix, some of which  were  from \cite{floudas1990collection,zheng2021t,globallib}, and the number of variables
        of the test problems ranges from 2 to 19 and the degree  from 2 to 40, as shown in Table \ref{tab:prob1-10}.

        For the experiments, the basic SOS relaxation method and the proposed block circulant SOS relaxation method were applied
         to obtain the SDP relaxation problem and the T-SDP relaxation problem, respectively.
        Then,  the T-SDP relaxation problem was transformed into an equivalent SDP problem as described in Section \ref{t-sdp}.
        For computation, we used Julia  1.7.3 with 
        Mosek \cite{aps2023mosek} 
        on a PC (Intel(R) Core(TM) i7-1185G7 @ 3.00GHz, 16GB, windows 10 Pro).

      Problems 1  \cite{floudas1990collection}, 2, 3, 4 \cite[Example 2]{zheng2021t} in Section \ref{result} demonstrate that the block circulant SOS 
      relaxation takes shorter computational time
        than the basic SOS relaxation method.
        Problem 5 \cite[st\_bpaf1b]{globallib} illustrates a case where the block circulant SOS relaxation problem is not feasible.
        With problem 6 \cite[Problem 2.9.1]{floudas1990collection}, we demostrate that the block circulant SOS relaxation can be an alternative
        approach to the basic SOS relaxation for the SDP relaxation problem with numerical instability, 
         which can be viewed as an additional benefit of the proposed T-SDP relaxation
         to its capability of handling larger-sized  problems described in Section~\ref{tsdp-relax}.
        For problem 7 \cite[st\_e34]{globallib}, we compare the
        basic SOS relaxation with the proposed block circulant SOS relaxation for large size problem.
        Problems 8, 9, 10 show the numerical efficiency of the proposed T-SDP relaxation.
        \begin{table}[htbp]
            \caption{Test problems.  $n$: the number of variable, $d$: the degree of $f$,
             $\max w_i$: the maximum of the degree of constraints, $N$: the relaxation level used to generate the SDP relaxation and the T-SDP relaxation,    $s(N) =  \binom{n+N}{N}$,
         and $s(N-\tilde{w}_i) = \binom{n+N-\tilde{w}_i}{N-\tilde{w}_i}$.}
            \label{tab:prob1-10}
            \begin{center}
                \small
                \begin{tabular}{lcccccc} \hline
                    No. & $n$ & $d$  &   $\max_{i} \{w_i\}$  & $N$  & $s(N)$ & $s(N-\tilde{w}_i)$  \\ \hline
                    1  \cite{floudas1990collection} & 10 & 2 & 1 & 2 &66 &11  \\
                    2 & 2 & 40 & 2 & 20 & $231$ & 210 \\
                    3 & 2 & 20 & 2 & 10 & $281$ & 220 \\
                     4 \cite[Example 2]{zheng2021t} & 2 & 58 & 2 & 29 & $465$ & 435 \\
                     5 \cite[st\_bpaf1b]{globallib} & 10 & 2 & 1 & 3 & $286$ & $66$ \\
                       6 \cite[Problem 2.9.1]{floudas1990collection} & 3 & 1 & 2 & 6 & $84$ & $56$ \\
                        7 \cite[st\_e34]{globallib} & 6 & 1 & 1 & 5 & $462$ & $210$ \\
                        8 & 11 & 6 & 2 & 3 & $364$ & $78$ \\
                        9 & 19 & 4 & 2 & 2 & $210$ & $20$ \\
                       10 & 19 & 4 & 2 & 2 & $210$ & $20$
                      \\ \hline
                \end{tabular}
            \end{center}
        \end{table}

        \subsection{Numerical results}
        \label{result}

        Tables~\ref{tab:prob1-and-5} and \ref{tab:prob6-and-10} report the numerical results
        on problems 1-5 and 6-10, respectively.
        In the tables,   ``Relax." denotes the SDP relaxation or T-SDP relaxation. In this subsection, the SDP relaxation 
        and the T-SDP relaxation correspond to 
        the basic SOS relaxation and the proposed block circulant SOS relaxation, respectively.
        ``Pn" (pattern) describes the 
         block sizes, more precisely, the sizes are arranged with the following columns ``$(m_0, l_0)$" and ``$(m_i, l_i)"$.
        For instance, $(m_0, l_0) = (66, 1)$ and \blue{ $(m_i, l_i) = (11, 1)$} indicates that $\phi_0 \in \Sigma^1[x]_{2N}$ and $\phi_i \in  \Sigma^1[x]_{2(N-\tilde{w}_i)}$ in the feasible set $\Gamma$ of the relaxation problem \eqref{p_k_tilde}
        which results in the T-SDP relaxation equivalent to the  basic SOS relaxation.
        For $(m_0, l_0) = (11, 6)$ and $(m_i, l_i) = (11, 1)$, we have
         $\phi_0 \in \Sigma^6[x]_{2N}$ and $\phi_i \in \Sigma^1[x]_{2(N-\tilde{w}_i)}$ in $\Gamma$.
         More precisely,  the $6$-block circulant SOS polynomial is employed to $\phi_0$ and an SOS polynomial to $\phi_i$.
         Tables~\ref{tab:prob1-and-5} and \ref{tab:prob6-and-10} also include the numbers and sizes
         of positive semidefinite matrices based on Table~\ref{tab:size-and-number}.
         ``\# of var." and ``\# of nnz" indicate
          the number of decision variables in the variable matrix and the number of nonzero elements of the SDP to be  solved, respectively, and 
          ``Opt.val" denotes its optimal value. ``CPU1", ``CPU2", and ``Tot." denote the time for generating the SDP problem to be solved, the computational time for solving
         the SDP problem, which are added to show the total time in seconds, respectively.

              \begin{table}[htbp]
                \centering
                \rotatebox{90}{
                    \begin{minipage}{\textheight}
                        \caption{Comparison of SDP  and T-SDP relaxation for Problems 1-5.}
                    \label{tab:prob1-and-5}
                                \begin{center}
                                    \footnotesize
                                    \begin{tabular}{ccccccccccc}

                            \hline
                            Relax. & Pn. & $(m_0, l_0)$ & $(m_i, l_i)$ &
                            No. * size of PSD matrix
                           &  \# of var. & \# of nnz & Opt.val & CPU1(s) & CPU2(s) & Tot.(s) \\ \hline
                            \hline
                            \multicolumn{11}{c}{Problem 1 (the known optimal value is 0.375)} \\
                            \hline
                            SDP & 1 & $(66,1)$ & $(11,1)$ &
                                [$1$ * $(66 \times 66)$,
                                $12$ * $(11 \times 11)$]
                            &
                            $3003$ & $8228$ & $0.37500$ & $2.99$ & $0.42$ & $3.41$ \\
                            T-SDP & 2 & $(11,6)$ & $(11,1)$ &
                            \begin{tabular}{lr}
                                [$4$ * $(22 \times 22)$,
                                $12$ * $(11 \times 11)$]
                            \end{tabular}
                            &
                            $1804$ & $45884$ & $0.37500$ & $2.22$ & $0.66$ & $2.88$ \\ \hline
                            \hline
                            \multicolumn{11}{c}{Problem 2 (the known optimal value is $14$)} \\
                            \hline
                            SDP & 1 & $(231,1)$ & $(210,1)$ & [$1$  * ($231\times 231$), $4$ * ($210 \times 210$)]  &$115416$ & $406161$ & $14.00000$ & $66.66$ & $13.14$ & $79.80$ \\
                            T-SDP & 2 & $(33,7)$ & $(210,1)$ &             [$4$ * $(66 \times 66)$,
                            $4$ * $(210 \times 210)$] &
                            $97464$ & $931200$ & $14.00000$ & $53.84$ & $15.89$ & $69.73$ \\
                            T-SDP & 3 & $(33,7)$ & $(105,2)$ &
                            [$4$ * $(66 \times 66)$,
                            $8$ * $105 \times 105$] &
                            $53364$ & $1269808$ & $14.00000$ & $32.61$ & $4.64$ & $37.25$ \\ \hline
                            \hline
                            \multicolumn{11}{c}{Problem 3 (the known optimal value is $1$)} \\
                            \hline
                            SDP & 1 & $(286,1)$ & $(220,1)$ &
                            [$1$*$(286 \times 286)$, $6$*$(220 \times 220)$] &
                            $186901$ & $662596$ & $1.00000$ & $187.55$ & $20.16$ & $207.71$ \\
                            T-SDP & 2 & $(143,2)$ & $(220,1)$ &
                            [$2$*$(143 \times 143)$, $6$*$(220 \time 220)$] &
                            $166452$ & $743413$ & $1.00000$ & $186.36$ & $19.30$ & $205.66$ \\
                            T-SDP & 3 & $(286,1)$ & $(110,2)$ &
                            [$1$*$(286 \times 286)$, $12$*($110 \times 110$)] &
                            $114301$ & $1229508$ & $1.00000$ & $110.60$ & $14.92$ & $125.52$ \\
                            T-SDP & 4 & $(143,2)$ & $(110,2)$ &
                            [$2$*$(143\times 143)$, $12$*$(110\times 110)$] &
                            $93852$ & $1310326$ & $1.00000$ & $99.61$ & $12.75$ & $112.36$ \\
                            T-SDP & 5 & $(26,11)$ & $(220,1)$ &
                            [$6$*$(52 \times 52)$, $6$*$(220 \times 220)$] &
                            $154128$ & $2082916$ & $1.00000$ & $176.65$ & $25.88$ & $202.53$ \\
                            T-SDP & 6 & $(26,11)$ & $(110,2)$ &
                            [$6$*$(52 \time 52)$, $12$*$(110 \times 110)$] &
                            $81528$ & $2649828$ & $1.00000$ & $88.60$ & $15.97$ & $104.57$ \\ \hline
                            \hline
                            \multicolumn{11}{c}{Problem 4 (the known optimal value is $1$)} \\
                            \hline
                            SDP & 1 & $(465,1)$ & $(435,1)$ &
                            [$1$*$(465\times 465)$, $2$*$(435\times 435)$] &
                            $298005$ & $1351575$ & $1.00000$ & $329.86$ & $64.16$ & $394.02$ \\
                            T-SDP & 2 & $(93,5)$ & $(435,1)$ &
                            [$3$*$(186 \times 186)$, $2$*$(435 \times 435)$] &
                            $241833$ & $2972864$ & $1.00000$ & $306.78$ & $63.17$ & $369.95$ \\
                            T-SDP & 3 & $(31,15)$ & $(435,1)$ &
                            [$8$*$(62 \times 62)$, $2$*$(435\times 435)$] &
                            $205284$ & $5633042$ & $1.00000$ & $221.96$ & $52.91$ & $274.87$ \\ \hline
                            \hline
                            \multicolumn{11}{c}{Problem 5 (the known optimal value is $-42.96256$)} \\
                            \hline
                            SDP & 1 & $(286,1)$ & $(66,1)$ &
                            [$1$*$(286 \times 286)$, $30$*$(66 \times 66)$] &
                            $107371$ & $438988$ & $-42.96256$ & $509.51$ & $421.66$ & $931.17$ \\
                            T-SDP & 2 & $(143,2)$ & $(66,1)$ &
                            [$2$*$(143 \times 143)$, $30$*$(66 \times 66)$] &
                            $86922$ & $520338$ & $-42.96255$ & $428.96$ & $449.61$ & $878.57$ \\
                            T-SDP & 3 & $(143,2)$ & $\left\{\begin{array}{c}(66,1) \\ (33,2) \end{array}\right.$ &
                            $\begin{array}{c}[ 2 * (143 \times 143), \\ 10*(66 \times 66), 40*(33 \times 33)] \end{array} $ &
                            $65142$ & $648402$ & $-42.96256$ & $334.09$ & $455.33$ & $782.42$ \\
                            T-SDP & 4 & $(143,2)$ & $(33,2)$ &
                            [$2$*$(143\times 143)$, $60$*$(33 \times 33)$] &
                            $54252$ & $869670$ & infeasible & - & - & - \\ \hline
                            \hline
                        \end{tabular}
                    \end{center}
                \end{minipage}
                } 
        \end{table}

        \begin{table}[htbp]
            \centering
            \rotatebox{90}{
                \begin{minipage}{\textheight}
                    \caption{Comparison of SDP  and T-SDP relaxation for Problems 6-10.}
                \label{tab:prob6-and-10}
                            \begin{center}
                                \footnotesize
                                \begin{tabular}{ccccccccccc}

                        \hline
                        Relax. & Pn. & $(m_0, l_0)$ & $(m_i, l_i)$ &
                        No. * size of PSD matrix
                       &  \# of var. & \# of nnz & Opt.val & CPU1(s) & CPU2(s) & Tot.(s) \\ \hline
                        \hline
                        \multicolumn{11}{c}{Problem 6 (the known optimal value is $-4$)} \\
                        \hline
                        SDP & 1 & $(84,1)$ & $(56,1)$ &
                      [$1$*$(84 \times 84)$, $8$*$(56 \times 56 $] &
                        $16338$ & $82320$ & $-3.99972$ & $4.30$ & $1.63$ & $5.93$ \\
                        T-SDP & 2 & $(42,2)$ & $(56,1)$ &
                        [$2$*$(42 \times 42)$, $8$*$(56\times 56)$] &
                        $14754$ & $89192$ & $-3.99998$ & $4.09$ & $1.45$ & $5.54$ \\
                        T-SDP &3 &  $(42,2)$ & $(28,2)$ &
                        [$2$*$(42 \times 42)$, $16$*$(28 \times 28)$] &
                        $8302$ & $160448$ & $-3.99999$ & $3.04$ & $1.47$ & $4.51$  \\ \hline
                        \hline
                        \multicolumn{11}{c}{Problem 7} \\
                        \hline
                        SDP & 1 & $(462,1)$ & $(210,1)$ &
                        [$1$*$(462 \times 462)$, $16$*$(210 \times 210)$] &
                        $461433$ & $2462544$ & $1.56195 \times 10^{-2}$ & $1877.58$ & $437.64$ & $2315.22$ \\
                        T-SDP & 2 & $(231,2)$ & $(210,1)$ &
                        [$2$*$(1231 \times 1231)$, $16$*$(210 \times 210)$] &
                        $408072$ & $2675042$ & $1.56195 \times 10^{-2}$ & $1225.32$ & $246.63$ & $1471.95$ \\
                        T-SDP &3 & $(231,2)$ &
                        $\left\{\begin{array}{c}(210,1) \\ (105,2) \end{array}\right.$ &
                        $\begin{array}{c} [2*(231 \times 231), \\
                            4*(210 \times 210), 24*(105 \times 105)]\end{array}$
                        & $275772$ & $4897672$ & $1.56195 \times 10^{-2}$ & $927.21$ & $305.72$ & $1232.93$ \\
                        T-SDP &4 &  $(42,11)$ &
                        $\left\{\begin{array}{c}(210,1) \\ (105,2) \end{array}\right.$ &
                        $\begin{array}{c} [6*(84 \times 84), \\
                            4*(210 \times 210), 24*(105 \times 105)]\end{array}$
                         & $243600$ & $8057606$ & $1.56195 \times 10^{-2}$ & $999.07$ & $387.34$ & $1386.41$ \\ \hline
                        \hline
                        \multicolumn{11}{c}{Problem 8} \\
                        \hline
                        SDP & 1 & $(364,1)$ & $(78,1)$ &
                        [$1$*$(364\times 364)$, $24$*$(78 \times 78)$] &
                        $140374$ & $546208$ & $-2.71766$ & $964.06$ & $300.58$ & $1264.64$ \\
                        T-SDP & 2 & $(182,2)$ & $(78,1)$ &
                        [$2$*$(91 \times 91)$, $24$*$(78 \times 78)$] &
                        $107250$ & $678036$ & $-2.71766$ & $542.06$ & $227.47$ & $769.53$ \\
                        T-SDP & 3 & $(182,2)$ &
                        $\left\{\begin{array}{c} (78,1) \\ (39,2) \end{array}\right.$ &
                        $\begin{array}{c}[2*(91 \times 91), \\
                            2*(78 \times 78), 44*(39 \times 39)]\end{array}$
                         & $73788$ & $928920$ & $-2.71766$ & $302.05$ & $331.75$ & $633.80$ \\
                         T-SDP & 4 & $(182,2)$ & $(39,2)$ &
                        [$2$*$(91 \times 91)$, $48$*$(39 \times 39)$] &
                        $70746$ & $1067368$ & $-2.71766$ & $374.34$ & $368.06$ & $742.40$ \\ \hline
                        \hline
                        \multicolumn{11}{c}{Problem 9} \\
                        \hline
                        SDP & 1 & $(210,1)$ & $(20,1)$ &
                        [$1$*$(210 \times 210)$, $40$*$(20\times 20)$] &
                        $30555$ & $90500$ & $-4.46273$ & $160.40$ & $46.52$ & $206.92$ \\
                        T-SDP & 2 & $(105,2)$ & $(20,1)$ &
                        [$2$*$(105\times 105)$, $40$*$(20 \times 20)$] &
                        $19530$ & $134062$ & $-8.27671$ & $64.45$ & $57.25$ & $121.70$ \\
                        T-SDP & 3 & $(105,2)$ & $(10,2)$ &
                        [$2$*$(105 \times 105)$, $80$*$(10 \times 10)$] &
                        $15530$ & $177874$ & infeasible & - & - & - \\
                        T-SDP & 4 & $(42,5)$ & $(20,1)$ &
                        [$3$*$(84 \times 84)$, $40$*$(20 \times 20)$] &
                        $19110$ & $439292$ & infeasible & - & - & - \\ \hline
                        \hline
                        \multicolumn{11}{c}{Problem 10} \\
                        \hline
                        SDP & 1 & $(210,1)$ & $(20,1)$ &
                        [$1$*$(210 \times 210)$, $40$*$(20 \times 20)$] &
                        $30555$ & $90500$ & $-1.00000$ & $91.22$ & $80.02$ & $171.24$ \\
                        T-SDP & 2 & $(105,2)$ & $(20,1)$ &
                        [$2$*$(105 \times 105)$, $40$*$(20 \times 20)$] &
                        $19530$ & $134062$ & $-1.00000$ & $78.67$ & $51.16$ & $129.83$ \\
                        T-SDP & 3 & $(42,5)$ & $(20,1)$ &
                        [$3$*$(84 \times 84)$, $40$*$(20 \times 20)$] &
                        $19110$ & $439292$ & $-1.00000$ & $71.22$ & $36.36$ & $107.58$ \\
                        T-SDP & 4 & $(21,10)$ & $(20,1)$ &
                        [$6$*$(42 \times 42)$, $40$*$(20 \times 20)$] &
                        $13818$ & $810308$ & $-1.00000$ & $52.20$ & $53.69$ & $105.89$ \\
                        T-SDP & 5 & $(42,5)$ & $(10,2)$ &
                        [$3$*$(84 \times 84)$, $80$*$(10 \times 10)$] &
                        $15110$ & $483104$ & $-1.00000$ & $52.68$ & $39.89 $ & $92.57$ \\ \hline
                    \end{tabular}
                \end{center}
            \end{minipage}
            } 
    \end{table}

    \subsubsection{Problem 1}  \label{Pr1result}

      In Table~\ref{tab:prob1-and-5}, we see that the number of decision variables for the proposed T-SDP relaxation method is $1804$ is smaller
      than $3003$ for the SDP relaxation method, resulting in shorter CPU1 (the
      computational time  to generate the SDP problem). 
      However,
      the T-SDP relaxation took slightly longer computation time in CPU2. This has been caused by the larger number of nonzero elements shown in the column
      of ``\# of nnz.''.
       Nevertheless, we can see in the ``Tot.'' column that the T-SDP relaxation consumed shorter overall computing time. 
        \subsubsection{Problem 2}


       For problem~2, Table~\ref{tab:prob1-and-5} displays the results by the SDP relaxation (Pn 1) and
       two T-SDP relaxations (Pn 2 and Pn 3). 
       The difference between Pn 2 and Pn 3  is that the block circulant SOS relaxation was applied
       only to $\phi_0$ for T-SDP (Pn 2)
       in the feasible set $\Gamma$ of \eqref{p_k_tilde} whereas the block circulant SOS relaxation was applied to both $\phi_0$ and $\phi_i$
       for  T-SDP (Pn 3).

        We observe in Table~\ref{tab:prob1-and-5} that the number of decision variables is significantly reduced to $97464$ and $53364$ in  T-SDP (Pn 2) and T-SDP (Pn 3),
        respectively,
        from $115416$ in the  SDP relaxation (Pn 1).
        As a result, the time for  generating the SDP problem to be solved for both T-SDP relaxations is reduced.
        However, the computational cost of solving the resulting T-SDP (Pn 2)  is slightly higher than the SDP relaxation for the same reason
         mentioned  in problem 1.
        On the other hand,
        the number of decision variables in  T-SDP (Pn 3)  is much smaller than that in T-SDP (Pn 2), and thus
        the computational time for solving T-SDP (Pn 3)  was  shorter than that for the SDP relaxation,
        even when the number of nonzero elements is taken into account.
        Hence, the total time for both T-SDP (Pn 2) and T-SDP (Pn 3) is less that that of the SDP relaxation (Pn 1).
        \subsubsection{Problem 3}


        The results for problem 3 show that 
        the number of variables in the T-SDP relaxations  is less than that of the SDP relaxation. Consequently, it took less time
        to solve the T-SDP relaxations than the SDP relaxation as shown in the last column.
       As  the number of variables in the T-SDP relaxations decreases,
       the number of nonzero elements increases from $662596$ in the SDP relaxation to $2649828$ for  T-SDP in the last row.
         The total times  for all T-SDP relaxations are shorter than  that for the SDP relaxation.
        \subsubsection{Problem 4}



        As shown in the column CPU1 on problem~4,  generating T-SDP relaxation problems consumed shorter computational time than
        the SDP relaxation. Also, two T-SDP relaxations, (Pn 2) and (Pn 3),  spent less computational time than that of 
        the SDP relaxation, despite the increase
        in the number of nonzero elements.

        \subsubsection{Problem 5}

	T-SDP relaxation problems may become infeasible depending on the values of $m_i$ and $l_i$, as illustrated   with this problem.
%
        In Table~\ref{tab:prob1-and-5},
          the  T-SDP relaxation (Pn 4)
           is obtained by using the $2$-block circulant SOS polynomial for $\phi_0$ and the $2$-block circulant SOS polynomial for $\phi_i$ which corresponds to the constraints for upper and lower bounds on each variable $x_1, \dots, x_{10}$,
          and SOS polynomials to $\phi_i$ which is associated  with the other inequality constraints in the feasible set $\Gamma$. 
          This is described as $(66,1)$ and $(33,2)$ 
          in the $(m_i,l_i)$ column.

  As discussed in Section~\ref{prop}, Assumption \ref{Assumption2} must hold for the T-SDP relaxation problem to have a feasible solution, but here we
  see that Assumption~\ref{Assumption2} no longer holds. 
        The T-SDP relaxation (Pn 4) becomes infeasible, as shown in the last row of Table~\ref{tab:prob1-and-5},
         when  the number of blocks $l_0$ or $l_i$ is increased in the block circulant SOS relaxation. 
      However,  
        the T-SDP relaxations (Pn 2) and (Pn 3)  are still feasible and computationally less expensive than the SDP relaxation.

        \subsubsection{Problem 6}

        We observe the effectiveness of  the  T-SDP  relaxation methods with this  problem.
         It is known that the optimal value  obtained by the SDP relaxation with relaxation level $N = 4$ is $-4$ \cite{lasserre2001global}. The optimal value 
          with $N = 6$, however, is $-3.99972$, which may contain numerical error.
          If the basic SOS relaxations with increasing the degree of the SOS polynomials is employed, then  the optimal value  by the SDP relaxation
          with a larger relaxation level approaches the optimal value of the original POP.
          We mention that the results in Table~\ref{tab:prob6-and-10} do not follow the theoretical result, which may be caused by the numerical instability of the SDP relaxation problem.
        In such cases, the T-SDP relaxation is still effective, and by changing parameters such as $l_0$ and $l_i$, another T-SDP relaxation problem
        can be generated.
        In particular,
        the smaller positive semidefinite matrices may improve the numerical accuracy.

        \subsubsection{Problem 7}


         The results  by the T-SDP relaxation (Pn 3) and (Pn 4) were obtained by applying the $2$-block and $11$-block circulant SOS polynomial, respectively. 
        The $2$-block circulant SOS polynomial for $\phi_i$
        corresponds to the constraints for upper and lower bounds on each variable $x_1, \dots, x_{6}$, and 
        the SOS polynomial for $\phi_i$ is associated
        with the other inequality constraints in the feasible set $K$ of \eqref{p_k_tilde}, as in 
        problem~5.

         In the rows for problem 7 in Table~\ref{tab:prob6-and-10}, we see  that
         the number of  variables in the
         T-SDP relaxation   is smaller than that of the SDP relaxation, which results in faster computational time for CPU1
         in all T-SDP cases.
        In particular, when the block circulant SOS polynomial was used to $\phi_i$ in the feasible set $\Gamma$,  a significant
         reduction in the number of decision variables and  CPU1 can be observed.
        For CPU2, as in the previous problems, there exists a trade-off between the increase in the number of nonzero elements and
        the decrease in the number of decision variables. Hence,
         increasing the parameters $l_0$ and $l_i$ does not necessarily reduce the computational time.
        We observe that all T-SDP relaxations outperform the SDP relaxation in terms of the total time.
        \subsubsection{Problem 8}

        We see  that the T-SDP relaxation  efficiently  provides an accurate optimal value for large-sized problems for this problem.
         The number of decision variables in the T-SDP relaxation is smaller that that of the SDP relaxation,
        which leads to shorter CPU1 time  in all cases. In particular,  CPU1 consumed by the
        T-SDP relaxation (Pn 3) is less than $\frac{1}{3}$ of that of the SDP relaxation.
        For CPU2, only the T-SDP relaxation (Pn 2) took less time   compared to the SDP relaxation, but 
        all T-SDP relaxations outperform the SDP relaxation in terms of the total time.

        \subsubsection{Problem 9}

	We discuss the performance of the T-SDP relaxation with two similar problems, 9 and 10.
	The objective functions of two problems are similar in that they
	 have the same  number of variables $n = 19$ and degree $d = 4$.
	The constraints of the
	problems are equivalent.  The T-SDP relaxation becomes infeasible for problem 9, while it provides an accurate optimal value
	for problem 10.


        In Table~\ref{tab:prob6-and-10}, the optimal value of the T-SDP relaxation (Pn 2) is smaller than that of the SDP relaxation, indicating that
        a weak lower bound was attained. In addition, the T-SDP relaxation 
        (Pn 3) and (Pn 4) were infeasible. Thus, the T-SDP relaxation may not work effectively for some problems, especially for the cases where the objective function does not satisfy  Assumption \ref{Assumption2}
         on the block circulant SOS polynomials.

        \subsubsection{Problem 10}

        We discuss cases where the T-SDP relaxation performs well for the
         problems with the same number of variables and degree as problem 9.
         The objective function of problem 10 contains a block circulant SOS polynomial.

        In the results for problem 10 in Table~\ref{tab:prob6-and-10}, we  see that the performance of the T-SDP relaxation 
        is better than the SDP relaxation,
        unlike problem 9. For CPU1,   all T-SDP relaxations spent shorter computational time than the  SDP relaxation. Furthermore, all T-SDP relaxations
        perform superior to 
        the SDP relaxation for CPU2.
  \subsection{Determining $l_i$} 
      \label{selection}
  From the numerical results in Section \ref{result},
we have observed that  the choices of $m_i$ and $l_i$ not only affect
   the computational time but also 
  the resulting objective function.
The sizes of the decomposed variable matrices  in Table~\ref{tab:size-and-number}, which are determined by  the values of $l_i$ among others,
 contribute to the computational time and can serve as estimates for it.

  Garksta et al.~\cite{Garstka2020} showed that 
  the computational time required to solve an SDP problem can be roughly estimated 
  by the cubes of the sizes of the decomposed matrices.
  For the $N$th-level relaxation, let $n_1, n_2, \dots, n_s$ be the sizes of PSD matrices 
  in Table~\ref{tab:size-and-number}.
 Specifically, with an approximation formula:
  \begin{align}
    c_0 + c_1 \sum_{j=1}^s n_j + c_2 \sum_{j=1}^s n_j^2 
    + c_3 \sum_{j=1}^s n_j^3,
    \label{regression-formula}
  \end{align}
  least squares can be employed to compute 
  the coefficients $c_0, c_1, c_2$ and $c_3$ 
  so that the approximation formula aligns with
  the computational time  observed in the numerical experiments.

\begin{figure}[tp]
      \begin{minipage}{0.48 \textwidth}
        \centering 
        \includegraphics[width=\textwidth]{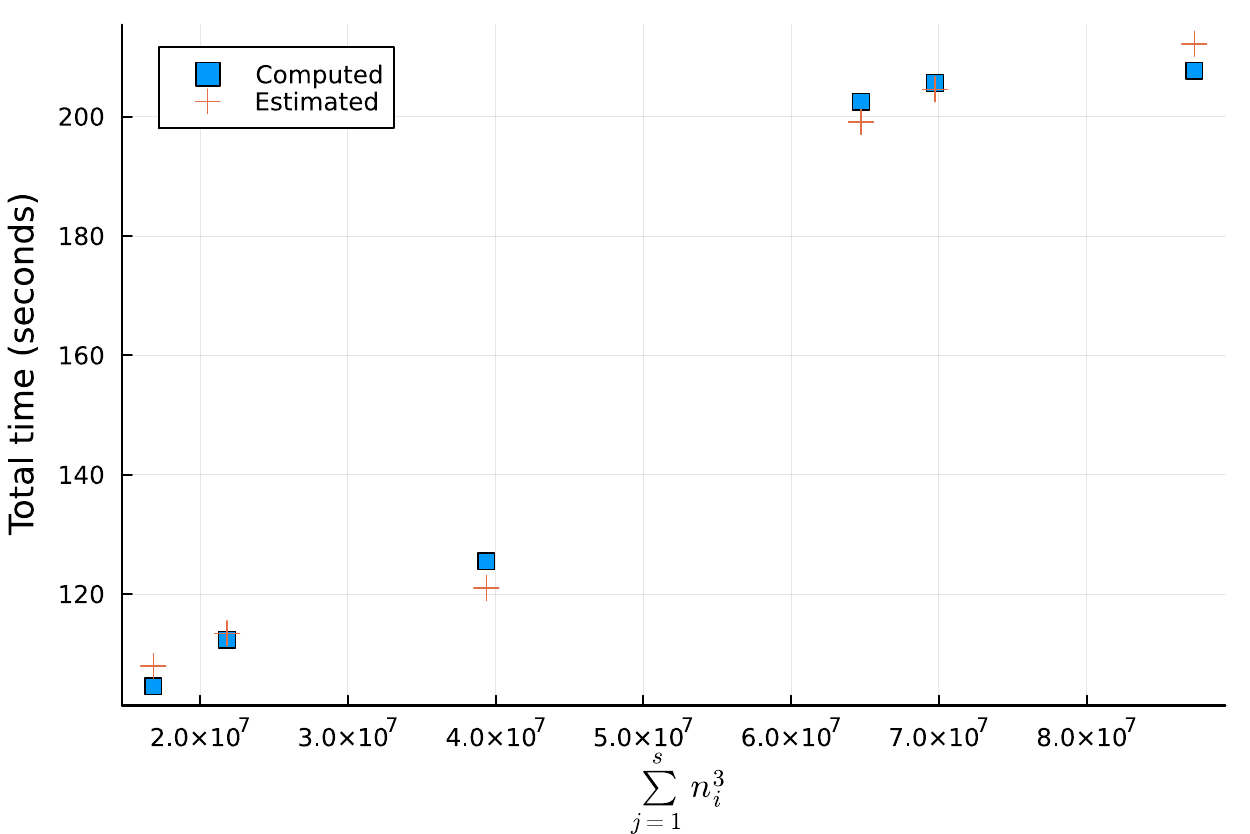}
        \newline
        Problem 3
      \end{minipage}
      \begin{minipage}{0.48 \textwidth}
        \centering
        \includegraphics[width=\textwidth]{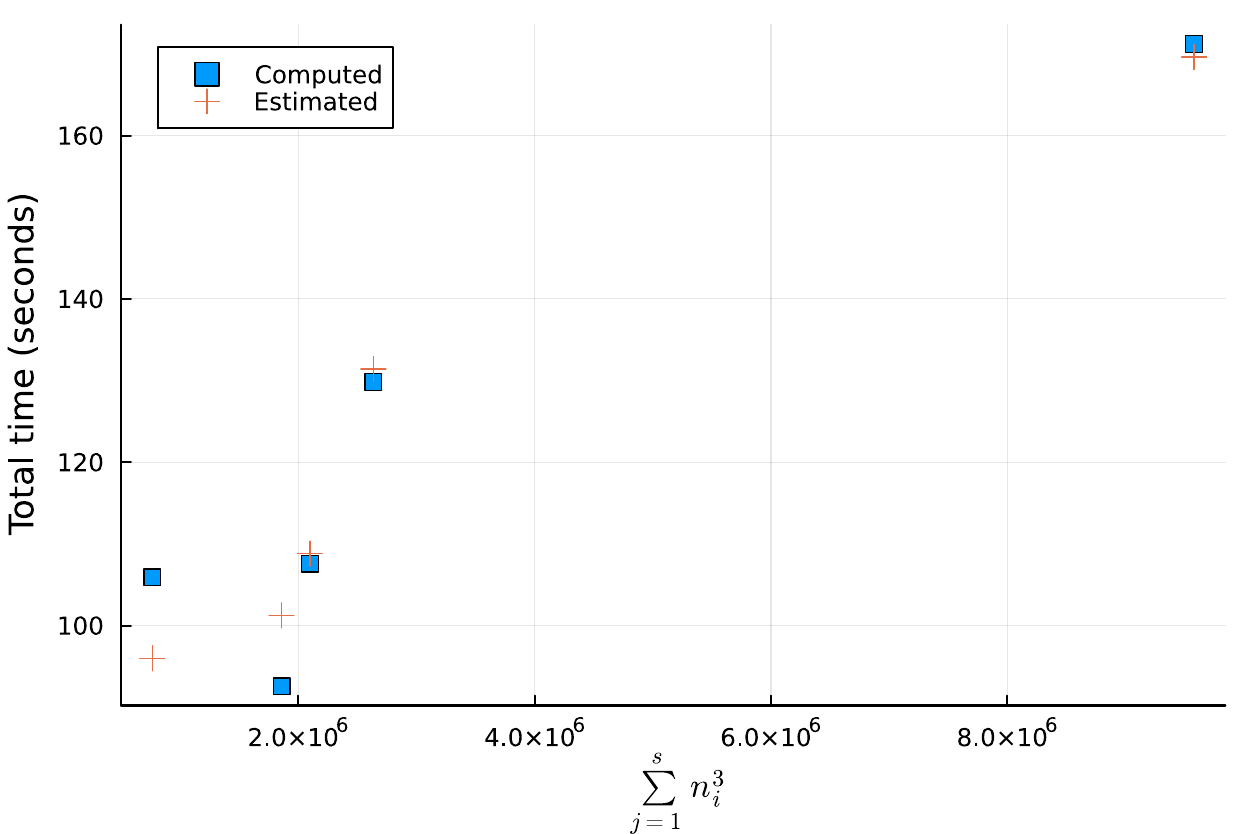}
        \newline
        Problem 10
      \end{minipage}
      
      \caption{Approximation for total computational time} \label{figure:fitting}
    \end{figure}

Figure~\ref{figure:fitting} compares the actual computational 
times for Problems 3 and 10 in 
Tables~\ref{tab:prob1-and-5} and \ref{tab:prob6-and-10}
and the estimated times from the approximation   
formula \eqref{regression-formula}.
We see that the formula \eqref{regression-formula} provides a
good approximation for the computational time  in the left figure of Figure~\ref{figure:fitting} for Problem~3. However, 
the accuracy of the approximation for Problem~10, displayed on the right figure of  Figure~\ref{figure:fitting},
 is comparatively lower.  
The values of $m_0, l_0, m_i$, and $l_i$  are constrained by the relations  \eqref{mlsplitting}, thus
restricting their freedom of choice, which in turn can affect the accuracy of approximation.
Nevertheless, Figure~\ref{figure:fitting} still shows a trend indicating that smaller values of 
$\sum_{j=1}^s n_s^3$
correspond to  shorter computational times.
Alternatively, the numerical results 
in Section \ref{result} 
suggest that 
decomposing the variable matrix into smaller-sized matrices  results in an inferior 
objective value. In addition, if the variable matrix is  decomposed 
into  very small-sized matrices, T-SDP can become infeasible 
as in the cases of Patterns 3 and 4 for Problem 9.
Hence,  a trade-off exists 
between the computational time and  objective value,
necessitating careful selection of $m_0, l_0, m_i$, and $l_i$.

%
    \section{Concluding remarks}
        \label{conc}
        We have proposed a relaxation method for constrained POPs using block circulant SOS polynomials. The proposed  block circulant SOS relaxation method is an extension of the SOS relaxation method proposed by Parrilo~\cite{parrilo2003semidefinite}, and
        also an extension of T-SDP relaxation  for
        unconstrained POPs  \cite{zheng2021t} in the sense that
	it  can be transformed into an equivalent T-SDP  using the third-order tensor.
	Theoretical analysis of the feasibility and global optimality of the T-SDP relaxation problem have been
	studied, and the convergence of the sequence of T-SDP relaxations with increasing the relaxation level has been
	established.

       The T-SDP relaxation method has shown to perform superior in computation to the SDP relaxation method from the basic SOS relaxation
       as it induces smaller-sized  problem.
         The computational results in Section~\ref{expe} show  that
         the T-SDP relaxation method provides the same quality of optimal values as the SDP relaxation method with shorter
         computational time under Assumption~\ref{Assumption2}.

        For future study, it will be interesting to  exploit sparsity for further improving the numerical efficiency.
        As in  \cite{waki2006sums},      the chordal sparsity of the block circulant SOS relaxation method can be developed for further reducing the computational cost.
        Another issue is to develop an SDP solver that directly solves T-SDP, instead of converting the T-SDP relaxation to an equivalent SDP relaxation.
        Then the computational time for the transformation can be saved.
       However, as mentioned earlier, the advantage of reducing the size of the relaxation problem
        when transforming from T-SDP to SDP may not be
       achieved in some case,
       which may require careful investigation.

\bibliographystyle{abbrv} 

\bibliography{ref.bib}

    \appendix
    \section*{Appendix}
{\scriptsize
        \textbf{Problem 1 \cite{floudas1990collection}}
        \begin{maxi*}
            {}  
            {\sum_{i=1}^{9} x_i x_{i+1} + \sum_{i=1}^{9} x_i x_{i+2} + x_1 x_7 + x_1 x_9 + x_1 x_{10} + x_2 x_{10} + x_4 x_7} 
            {}  
            {}  
            \addConstraint{\sum_{i=1}^{10} x_i}{=1} 
            \addConstraint{x_i}{\ge 0,}{\quad i=1, \dots, 10.}
        \end{maxi*}

        \textbf{Problem 2}
        \begin{mini*}{}{x_1^{40} + x_1^{18}x_2^{20} + x_1^{14}x_2^{24} + x_1^{26}x_2^{8} + x_1^{22}x_2^{12} + x_1^{26}x_2^{4}}{}{}
        \breakObjective{+ x_1^{22}x_2^{8} + x_1^{18}x_2^{8} + x_1^{14}x_2^{12} + x_1^{2}x_2^{20} + x_1^{20} + x_1^{12}x_2^{2} + x_1^{8}x_2^{6} + x_2^{2}}
        \addConstraint{x_1, x_2}{\in \{-1, +1\}}
        \end{mini*}

        \textbf{Problem 3}
        \begin{mini*}
            {}  
            {x_1^{20} + x_2^{2}x_3^{2} - 2.0x_2x_3^{3} + x_3^{4} - 4.0x_2x_3^{2} + 4.0x_3^{3} + 4.0x_3^{2}} 
            {}  
            {}  
            \addConstraint{x_1, x_2, x_3}{\in \{-1, +1\}} 
        \end{mini*}

        \textbf{Problem 4 \cite[Example 2]{zheng2021t}}
        \begin{maxi*}
            {}{x_1^{10}x_2^4 + x_1^8x_2^{12} + x_1^{24}x_2^2 + x_1^{24}x_2^6 + x_1^{32}x_2^2 + x_1^8x_2^{28} + x_1^{28}x_2^{12}}{}{}
            \breakObjective{+ x_1^{10}x_2^{32} + x_1^{42}x_2^4 + x_1^{30}x_2^{18} + x_1^{20}x_2^{30} + x_1^{12}x_2^{40} + x_1^6x_2^{48} + x_1^2x_2^{54} + x_2^{58}}
            \addConstraint{x_1^2 + x_2^2}{= 1.}
        \end{maxi*}

        \textbf{Problem 5 \cite[st\_bpaf1b]{globallib}}
        \begin{mini*}
            {}  
            {x_1 x_6 + 2 x_1 - 2 x_6 + x_2 x_7 + 4 x_2 - x_7 + x_3 x_8} 
            {}  
            {}  
            \breakObjective{+ 8 x_3 - 2 x_8 + x_4 x_9 - x_4 - 4 x_9 + x_5 x_{10} - 3 x_5 + 5 x_{10}}
            \addConstraint{-8x_1 - 6x_3 + 7x_4 - 7x_5}{\le 1}
            \addConstraint{-6x_1 + 2x_2 + 3x_3 - 9x_4 - 3x_5}{\le 3}
            \addConstraint{6x_1 - 7x_3 - 8x_4 + 2x_5}{\le 5}
            \addConstraint{-x_1 + x_2 - 8x_3 - 5x_5}{\le 4}
            \addConstraint{4x_1 - 7x_2 + 4x_3 + 5x_4 + x_5}{\le 0}
            \addConstraint{5x_7 - 4x_8 + 9x_9 - 7x_{10}}{\le 0}
            \addConstraint{7x_6 + 4x_7 + 3x_8 + 7x_9 + 5x_{10}}{\le 7}
            \addConstraint{6x_6 + x_7 - 8x_8 + 8x_9}{\le 3}
            \addConstraint{-3x_6 + 2x_7 + 7x_8 + x_{10}}{\le 6}
            \addConstraint{-2x_6 - 3x_7 + 8x_8 + 5x_9 - 2x_{10}}{\le 2}
            \addConstraint{0 \le x_i}{\le 20, \quad i=1, \dots, 10.}
        \end{mini*}

        \textbf{Problem 6 \cite[Problem 2.9.1]{floudas1990collection}}
        \begin{mini*}{}{-2x_1 + x_2 - x_3}{}{}
        \addConstraint{x^\mathsf{T}B^\mathsf{T}Bx - 2r^\mathsf{T}Bx + ||r||^2 - 0.25||b - v||^2}{\ge 0}
        \addConstraint{x_1 + x_2 + x_3}{\le 4}
        \addConstraint{x_1}{\le 2}
        \addConstraint{x_3}{\le 3}
        \addConstraint{3x_2 + x_3}{\le 6}
        \addConstraint{x_i}{\ge 0, \quad i=1, \dots, 3.}
        \end{mini*}

        Here, 
        \begin{equation*}
            B =
            \left[ \begin{array}{@{\,}ccc@{\,}}
                 0 &  0 &  1 \\
                 0 & -1 &  0 \\
                -2 &  1 & -1
            \end{array} \right], \quad
            b =
            \left[ \begin{array}{@{\,}c@{\,}}
                3 \\
                0 \\
                4
            \end{array} \right], \quad
            v =
            \left[ \begin{array}{@{\,}c@{\,}}
                0 \\
                -1 \\
                -6
            \end{array} \right], \quad
            r =
            \left[ \begin{array}{@{\,}c@{\,}}
                1.5 \\
                -0.5 \\
                -5
            \end{array} \right].
        \end{equation*}

        \textbf{Problem 7 \cite[st\_e34]{globallib}}
        \begin{equation*}
            \begin{aligned}
                &\textrm{minimize} & 4.3x_1 + 31.8x_2 + 63.3x_3 + 15.8x_4 + 68.5x_5 + 4.7x_6 \\
                &\textrm{subject to} &
                17.1x_1 - 169x_1x_3 + 204.2x_3 - 3580x_3x_5 + 623.4x_5 - 3810x_4x_6 \\&& + 212.3x_4 + 1495.5x_6 - 18500x_4x_6 + 38.2x_2 \ge 4.97, \\
                & & 17.9x_1 - 139x_1x_3 + 113.9x_3 - 2450x_4x_5 + 169.7x_4 + 337.8x_5 \\ && - 16600x_4x_6 + 1385.2x_6 - 17200x_5x_6 + 36.8x_2 \ge -1.88, 
                  \end{aligned}
                 \end{equation*}   
                     \begin{equation*}
                    \begin{aligned}
                & & 26000x_4x_5 - 70x_4 - 819x_5 - 273x_2 \ge -69.08, \\
                & & 159.9x_1 - 14000x_1x_6 + 2198x_6 - 311x_2 + 587x_4 + 391x_5 \ge -118.02, \\
                & & 0 \le x_1 \le 0.31, \\
                & & 0 \le x_2 \le 0.046, \\
                & & 0 \le x_3 \le 0.068, \\
                & & 0 \le x_4 \le 0.042, \\
                & & 0 \le x_5 \le 0.028, \\
                & & 0 \le x_6 \le 0.0134. \\
            \end{aligned}
        \end{equation*}

        \textbf{Problem 8}
        \begin{mini*}
            {x \in \Real^{11}}  
            {\sum_{0 < \alpha \le 6} \mathrm{rand}(0:1)_\alpha x^\alpha} 
            {}  
            {}
            \addConstraint{\sum_{i=1}^{11} x_i^2}{= 1}
            \addConstraint{-1 \le x_i}{\le 1, \quad i=1, \dots, 11,}
        \end{mini*}
        where $\mathrm{rand}(0:1)_\alpha$ is a random number of $0$ or $1$.

        \textbf{Problem 9}
        \begin{mini*}
            {x \in \Real^{19}}  
            {\sum_{0 < \alpha \le 4} \mathrm{rand}(0:1)_\alpha x^\alpha} 
            {}  
            {}
            \addConstraint{\sum_{i=1}^{19} x_i^2}{= 1}
            \addConstraint{-1 \le x_i}{\le 1, \quad i=1, \dots, 19,}
        \end{mini*}
        where $\mathrm{rand}(0:1)_\alpha$ is a random number of $0$ or $1$ as in Problem 8.

        \textbf{Problem 10}
        \begin{mini*}
            {}  
            {x_1^4 + x_1^2x_{12}^2 + 2.0x_1^2x_{12}x_{17} + x_1^2x_{17}^2 + 2.0x_1x_2x_3x_{12} + 2.0x_1x_2x_3x_{17} + x_2^2x_3^2}
            {}  
            {}
            \breakObjective{+ x_2^2x_4^2 + x_2^2x_{15}^2 + 2.0x_2x_3^2x_{15} + 2.0x_2x_3x_7x_{15} + x_3^4 + 2.0x_3³x_7 + x_3^2x_7^2 + x_3^2x_8^2}
            \breakObjective{+ x_3^2x_{19}^2 + 2.0x_3x_4x_8x_{19} + 2.0x_3x_4x_{12}x_{19} + x_4^2x_8^2 + 2.0x_4^2x_8x_{12} + x_4^2x_{12}^2 + x_4^2x_{13}^2}
            \breakObjective{+ x_5^2x_9^2 + 2.0x_5^2x_9x_{14} + 2.0x_5^2x_9x_{18} + x_5^2x_{14}^2 + 2.0x_5^2x_{14}x_{18} + x_5^2x_{18}^2 + x_5^2x_{19}^2}
            \breakObjective{+ x_6^2x_{16}^2 + 2.0x_6x_7x_8x_{16} + 2.0x_6x_7x_{12}x_{16} + x_7^2x_8^2 + 2.0x_7^2x_8x_{12} + x_7^2x_{12}^2 + x_7^2x_{13}^2}
            \breakObjective{+ x_8^2x_{12}^2 + 2.0x_8^2x_{12}x_{17} + x_8^2x_{17}^2 + 2.0x_8x_9x_{10}x_{12} + 2.0x_8x_9x_{10}x_{17} + x_9^2x_{10}^2}
            \breakObjective{+ x_9^2x_{11}^2 + x_{10}^2x_{12}^2 + 2.0x_{10}^2x_{12}x_{17} + x_{10}^2x_{17}^2 + 2.0x_{10}x_{11}x_{12}^2 + 2.0x_{10}x_{11}x_{12}x_{17} }
            \breakObjective{+ x_{11}^2x_{12}^2 + x_{11}^2x_{13}^2 + x_{12}^2x_{16}^2 + 2.0x_{12}x_{13}x_{14}x_{16} + 2.0x_{12}x_{13}x_{16}x_{18} + x_{13}^2x_{14}^2}
            \breakObjective{+ 2.0x_{13}^2x_{14}x_{18}+ x_{13}^2x_{18}^2 + x_{13}^2x_{19}^2 + x_{15}^2x_{19}^2 + 2.0x_{15}x_{17}^2x_{19} + 2.0x_{15}x_{18}x_{19}^2}
            \breakObjective{+ x_{17}^4 + 2.0x_{17}^2x_{18}x_{19} + x_{18}^2x_{19}^2 + x_{19}^4 + x_{11}^2 + 2.0x_{11}x_{16} + x_{16}^2 + 2.0x_{11} + 2.0x_{16}}
            \addConstraint{\sum_{i=1}^{19} x_i^2}{= 1}
            \addConstraint{-1 \le x_i}{\le 1, \quad i=1, \dots, 19,}
        \end{mini*}
  }

\end{document}